\documentclass[11pt]{amsart}
\usepackage{fullpage, amsmath , amsfonts, amscd, epsfig, algorithmic, algorithm, amssymb, pstricks,pst-node, pstcol}

\psset{unit=1pt,arrowsize=4pt,linewidth=.7pt}
\psset{linecolor=blue}
\newgray{grayish}{.90}
\newrgbcolor{emgreen}{0 .5 0}
\def\vblack(#1,#2)#3{\cnode*[linecolor=black](#1,#2){3}{#3}}
\def\vwhite(#1,#2)#3{\cnode[linecolor=black,fillcolor=white,fillstyle=solid](#1,#2){3}{#3}}

\theoremstyle{plain}
  \newtheorem{theorem}{Theorem}
  \newtheorem{proposition}[theorem]{Proposition}
  \newtheorem{lemma}[theorem]{Lemma}
  \newtheorem{corollary}[theorem]{Corollary}
  \newtheorem{conjecture}[theorem]{Conjecture}
\theoremstyle{definition}
  \newtheorem{definition}[theorem]{Definition}
  
\theoremstyle{remark}
  \newtheorem{remark}[theorem]{Remark}
\theoremstyle{Problem}

\newcommand{\field}[1]{\mathbb{#1}}
\newcommand{\LE}{%
 \hbox{%
  \vbox{\hrule width 0.35em height 0.04 em}%
  \vbox{\offinterlineskip%
   \hbox{\kern -0.02em\vrule height 0.65em width 0.04em \hspace{0.1em}}%
  }%
 }%
}

\def\A{\mathcal{A}}
\def\B{\mathcal{B}}
\def\I{\mathcal{I}}

\def\F{\mathcal{F}}
\def\M{\mathcal{M}}

\def\R{\field{R}}
\def\Q{\field{Q}}
\def\S{\mathcal{S}}

\title{$\LE$-diagrams and totally positive bases inside the nonnegative Grassmannian}
\author{Suho Oh}
\address{Department of Mathematics, Massachusetts Institute of Technology, 
        77 Massachusetts Ave, Cambridge, MA 02139}      

\date{August 30, 2008}
\thanks{The author was supported in part by Samsung Scholarship.}

\begin{document}

\begin{abstract}
There is a cell decomposition of the nonnegative Grassmannian. For each cell, totally positive bases(TP-bases) is defined as the minimal set of Pl\"ucker variables such that all other nonzero Pl\"ucker variables in the cell can be expressed in those variables in a subtraction-free rational function. This is the generalization of the TP-bases defined for nonnegative part of $GL_k$ defined in \cite{FZ5}.

For each cell, we have a $\LE$-diagram and a natural way to label the dots inside the diagram with Pl\"ucker variables. Those set of Pl\"ucker variables form a TP-bases of the cell. Using mutations coming from 3-term Pl\"ucker relation, we conjecture that they can be mutated to a special set of Pl\"ucker variable $\S$. All other nonzero Pl\"ucker variables in the cell will be expressed as a subtraction-free Laurent polynomial in variables of $\S$. We define TP-diagrams to express the transformation procedure in terms of moves on a diagram.

We will prove the conjecture for certain class of cells called weakly-connected cells. Then we will study the connection with cluster algebras through lattice-path-matroid cells.

\end{abstract}

\maketitle

\section{Introduction}

Cluster algebras were first defined by Fomin and Zelevinsky to study total positivity in algebraic groups and canonical bases in quantum groups. In a naive sense, a cluster algebra is a ring generated by a set of generators called cluster variables. The cluster variables are grouped into seeds. The cluster variables are obtained recursively from initial seeds by a process called mutation. One of the nice properties of cluster algebras is that given any seed, all other cluster variables can be expressed as a Laurent polynomial in variables of the seed. The nonegativity conjecture states that the Laurent polynomial is subtraction-free. We will study a related problem.

Fomin-Zelevinsky (\cite{FZ5}) investigated the decomposition of the totally nonnegative part of $GL_k$ into cells, called the double Bruhat cells. They defined totally positive bases(TP-bases) as the set of minimal set of minors $\B$ such that all other minors can be expressed as subtraction-free rational function in terms of minors in $\B$. They also showed that minors coming from double-wiring diagrams form a TP-bases. It is conjectured in \cite{FZ6} that in each cell, all minors can be expressed as subtraction-free Laurent polynomial in any set of minors coming from any configuration of a double-wiring diagram of the cell. We want to generalize this.

The totally nonnegative Grassmannian $Gr_{kn}^{tnn}$ is defined as the set of elements in the Grassmannian $Gr_{kn}$ with all nonnegative Pl\"ucker coordinates. This was motivated from the classical notion of totally positive matrix, which is a matrix with all positive minors. It is well known that there is a decomposition of $Gr_{kn}$ into cells $S_{\M}$, the matroid strata. We can take the intersection of $S_{\M}$ with $Gr_{kn}^{tnn}$ to obtain a cell decomposition. Denote them $S_{\M}^{tnn}$ and we call the matroids such that $S_{\M}$ is nonempty as positroids. We call the cells as positroid cells.
 
 For each cells, we will study the TP-bases. TP-bases for each cell is defined as the set of algebraically independent Pl\"ucker variables $\B$ such that all other nonzero variables in the cell can be expressed as subtraction-free rational function in variables of $\B$.

 Be aware that TP-bases we use are different from TP-basis in \cite{DKK}. In \cite{DKK}, TP stands for tropical Pl\"ucker, while our usage of TP stands for totally positive. Nonetheless, they are related, and we show how in the last section.

For each cell, there is a combinatorial object called $\LE$-diagram. There is a canonical way to label the dots inside the $\LE$-diagrams with Pl\"ucker variables. Let $\S$ be the set of such Pl\"ucker variables, then $\S$ is a TP-bases. We call them the canonical TP-bases of the cell. Motivated from seed mutations of the cluster algebra, we will use Pl\"ucker mutation induced from 3-term-Pl\"ucker relations on our TP-bases. We say $\S$ and $\S'$ are mutation equivalent if $\S'$ can be obtained from $\S$ by Pl\"ucker mutations.

For a $\LE$-diagram, the boundary forms a lattice path from the upper-right corner to lower-left corner. Label the $n$ steps in this path by numbers $1,\cdots,n$ consecutively. Define $I=I(\lambda)$ as the set of lables of $k$ vertical steps in the path. Put dots on each edge of the boundary path. Connect all dots on same row and connect all dots on the same column. Then we get a $\LE$-graph. On a $\LE$-graph, direct all horizontal edges leftward and vertical edges downward. Put a weight on each edge, but let the value be $1$ for all vertical edges. Denote weight of each horizontal edge $e$ by $x_e$. We call this a $\LE$-network. Then this $\LE$-network completely parametrizes the cell(\cite{P}).

 Let $I = \{i_1,\cdots,i_k\}, [n]\setminus I = \{j_1,\cdots,j_{n-k}\}$ such that $i_1 < \cdots < i_k, j_1 > \cdots > j_{n-k}$. The Pl\"ucker variable $\Delta_{X,Y}$ denotes $$\Delta_{ (I \setminus \{ i_{x_1},\cdots,i_{x_k} \}) \cup \{j_{y_1},\cdots,j_{y_k} \}}$$ where $X = \{x_1,\cdots,x_r\},Y=\{y_1,\cdots,y_r\}, r \leq k$. Setting $\Delta_{\phi,\phi}=\Delta_{I(\lambda)}=1$, we get $$\Delta_{H,J} = \sum_{P: H \rightarrow J} \prod_{e \in P} x_e$$ where $P: H \rightarrow J$ indicates non-vertice-crossing paths from source $H$ to sinks $J$.

We conjecture the following.

\begin{conjecture}
\label{con:main}
For any cell in the nonnegative Grassmannian, canonical TP-bases $\S$ is mutation equivalent to a TP-bases $\S'$ such that all edge weights of $\LE$-network can be expressed as Laurent monomial in variables of $\S'$. This implies that all nonzero Pl\"ucker variables in the cell can be expressed as a subtraction-free Laurent polynomial in variables of $\S'$.
\end{conjecture}

 This will be proved for certain class of cells called weakly-connected cells. We will construct an algorithm to go from $\S$ to $\S'$ and every Pl\"ucker mutation we use will be expressed as a move on the $\LE$-diagram.

\begin{theorem}
For any weakly-connected cell in the nonnegative Grassmannian, canonical TP-bases $\S$ is mutation equivalent to a TP-bases $\S'$ such that all edge weights of $\LE$-network can be expressed as Laurent monomial in variables of $\S'$. This implies that all nonzero Pl\"ucker variables in the cell can be expressed as a subtraction-free Laurent polynomial in variables of $\S'$.
\end{theorem}

Now let's see how this is related to cluster algebras. The seeds and the mutation process can be expressed in terms of quivers. For each cell, we can get a quiver from the $\LE$-diagram. And from that quiver, we can assign a cluster algebra. Call this the canonical cluster algebra of the cell. We have the following result.

Pick a skew Young diagram inside $k$-by-$(n-k)$ rectangle. Denote the shape by $\lambda / \mu$. For Young diagram $\lambda$, put dots inside boxes in $\lambda / \mu$. Then we get a $\LE$-diagram. There is a cell associated to the $\LE$-diagram. This cell is a special case of set of weakly-connected cells. Such cells come from lattice-path-matroids. Let's say a box in the diagram is at $(x,y)$ if it is at $x$-th row and $y$-th column.

\begin{theorem}
For each box $(a,b)$, denote $k_{a,b}$ as the maximal number such that there is a box inside the diagram at $(a+k_{a,b},b+k_{a,b})$. Define $X_{a,b}$ as $\{a,a+1,\cdots,a+k_{a,b} \}$ and $Y_{a,b}$ as $\{b, b+1, \cdots, b+k_{a,b} \}$.

Let $\S'$ be set of Pl\"ucker variables of form $\Delta_{X_{a,b},Y_{a,b}}$ for all dots in the $\LE$-diagram. Set $x_e$ where $e$ is an edge heading towards $(x,y)$ as $\Delta_{X_x,Y_y}\Delta_{X_{x+1},Y_{y+2}} / \Delta_{X_x,Y_{y+1}}\Delta_{X_{x+1},Y_{y+1}}$. $\Delta_{H,J}$ where $H \not \subseteq [n]$ or $J \not \subseteq [n]$ is considered to be $1$.

$\S'$ is a TP-basis that is mutation equivalent to the canonical TP-basis of the cell. Furthermore, $\S'$ is a cluster inside the cluster algebra associated to the cell.

\end{theorem}

How we prove the theorem is the following. Since the cells we are looking at are weakly connected, our algorithm for TP-bases gives us $\S'$ from the canonical TP-bases. And we will show that for such cells described in the theorem, our transformation process is exactly the seed mutation of the associated cluster algebra. To show this, we use the transformation process for plabic graphs described in \cite{P}.

The paper is organized as follows. In section 2, we go over the basics of nonnegative Grassmannian and $\LE$-diagrams. In section 3, we go over the basics of cluster algebras.
In section 4, we go over the basics of plabic graphs. In section 5, we go over the basics of $\LE$-networks and how they parametrize the cells of nonnegative Grassmannian. In section 6, we define TP-bases and state the main theorem. We also state an algorithm such that when applied recursively, proves the theorem. In section 7, we show an example of how the algorithm works. In section 8,9 and 10 is devoted to proving the main lemma. In section 11, we prove the main theorem and show how to obtain $\S'$ from $\S$ using the algorithm. In section 12, we apply our algorithm on lattice-path-matroid cells and connect TP-bases to cluster algebras. In section 13, we provide some open problems.

\medskip

\textbf{Acknowledgment} I would like to thank my adviser, Alexander Postnikov for introducing me to the field, and helping me with useful advices. 

\section{Non-negative Grassmannian, positroids and $\LE$-diagrams}

Ordering $<_w, w \in S_n$ is defined as $a<_w b$ if $w^{-1} a < w^{-1} b$ for $a,b \in [n]$.

\begin{definition}
Let $A,B \in \ {[n] \choose k}$, $w \in S_n$ where 
$$ A = \{i_1, \cdots, i_k\}, i_1 <_w i_2 <_w \cdots <_w i_k $$
$$ B = \{j_1, \cdots, j_k\}, j_1 <_w j_2 <_w \cdots <_w j_k $$
Then we set $A \leq_w B$ if and only if $i_1 \leq_w j_1, \cdots, i_k \leq_w j_k$. This is called the Gale ordering on ${[n] \choose k}$ induced by $w$. We denote $\leq_t$ for $t \in [n]$ as $<_{c^{t-1}}$ where $c=(1,\cdots,n) \in S_n$.
\end{definition}

We can define matroids from above ordering, see \cite{G},\cite{BGW}. 

\begin{definition}
\label{def:mat2}
Let $\M \subseteq {[n] \choose k}$. Then $\M$ is a matroid if and only if $\M$ satisfies the following property. For every $w \in S_n$, the collection $\M$ contains a unique member $A \in \M$ maximal in $\M$ with respect to the partial order $\leq_w$.
\end{definition}
 
Now we can define a \textit{Schubert matroid} using the partial order $\leq_w$.

\begin{definition}
 For $I=(i_1, \cdots, i_k)$, the \textit{Schubert Matroid} $SM^{w}_I$ consists of bases $H=(j_1, \cdots, j_k)$ such that $I \leq_w H$.
\end{definition}

Now let's recall the Grassmannian. For detailed explanation, see \cite{F}. An element in the Grassmannian $Gr_{kn}$ can be understood as a collection of $n$ vectors $v_1, \cdots, v_n \in \R^k$ spanning the space $\R^k$ modulo the simultaneous action of $GL_k$ on the vectors. The vectors $v_i$ are the columns of a $k \times n$-matrix $A$ that represents the element of the Grassmannian. Then an element $V \in Gr_{kn}$ represented by $A$ gives the matroid $\M_V$ whose bases are the $k$-subsets $I \subset [n]$ such that $\Delta_I (A) \not = 0$.

Then $Gr_{kn}$ has a subdivision into \textit{matroid strata} $S_{\M}$ labeled by some matroids $\M$:
$$ S_{\M} := \{V \in Gr_{kn} | \M_{V} = \M \}$$
The elements of the stratum $S_{\M}$ are represented by matrices $A$ such that $\Delta_I (A) \not = 0$ if and only if $I \in \M$. Now we define the Schubert matroids, which corresponds to the cells of the matroid strata. 

Let us define the totally nonnegative Grassmannian, its cells and positroids.

\begin{definition}[\cite{P}] The \textit{totally nonnegative Grassmannian} $Gr^{tnn}_{kn} \subset Gr_{kn}$ is the quotient $Gr^{tnn}_{kn} = GL^{+}_{k} \backslash Mat^{tnn}_{kn}$, where $Mat^{knn}_{kn}$ is the set of real $k\times n$-matrices $A$ of rank $k$ with nonnegative \textit{maximal} minors $\Delta_I (A) \geq 0$ and $GL^{+}_k$ is the group of $k \times k$-matrices with positive determinant.
\end{definition}

\begin{definition}[\cite{P}] The \textit{totally nonnegative Grassmann cells} $S^{tnn}_{\M}$ in $Gr^{tnn}_{kn}$ is defined as $S^{tnn}_{\M} := S_{\M} \cap Gr^{tnn}_{kn}$. $\M$ is called a \textit{positroid} if the cell $S^{tnn}_{\M}$ is nonempty.
\end{definition}

Note that from above definitions, we get
$$ S^{tnn}_{\M} = \{ GL^{+}_k \bullet A \in Gr^{tnn}_{kn} | \Delta_I (A) >0 \textbf{ for } I \in \M, \Delta_I (A) = 0 \textbf{ for } I \not \in \M  \} $$

\begin{remark}
\label{rem:dbruhatcell}
For $k$-by-$n$ matrix $A$ such that $A_{\{1,\cdots,k\}} = ID_k$, define $\phi(A)=B$ where $B=(b_{ij})$ is a $k$-by-$(n-k)$-matrix such that $b_{ij} = (-1)^{k-j} a_{i+k,j}$. Then $\phi$ induces isomorphism between totally nonnegative part of $Gr_{k,2k}$ such that $\Delta_{\{1,\cdots,k\}} > 0$ and $\Delta_{\{k+1,\cdots,2k\}} > 0$ and the totally nonnegative part of $GL_k$. Moreover, it gives isomorphisms between double Bruhat cells in $GL_k$ and cells $S_{\M}^{tnn} \subset Gr_{k,2k}^{tnn}$ such that $\{1,\cdots,k\},\{k+1,\cdots,2k\} \in \M$.
\end{remark}

In \cite{P}, Postnikov showed a bijection between each positroid cells and an combinatorial object called Grassmann necklace. He also showed that those necklaces can be represented as objects called decorated permutations. Let's first see how they are defined.

\begin{definition}[\cite{P}]
A \textit{Grassmann necklace} is a sequence $\I = (I_1, \cdots, I_n)$ of subsets $I_r \subseteq [n]$ such that, for $i \in [n]$, if $i \in I_i$ then $I_{i+1}= (I_i \setminus \{i \}) \cup \{j\}$, for some $j \in [n]$; and if $i \in I_i$ then $I_{i+1} = I_i$. (Here the indices are taken modulo $n$.) In particular, we have $|I_1| = \cdots = |I_n|$.

\end{definition}

\begin{definition}[\cite{P}]
 A decorated permutation $\pi^{:} = (\pi, col)$ is a permutation $\pi \in S_n$ together with a coloring function $col$ from the set of fixed points $\{i | \pi(i) = i\}$ to $\{1,-1\}$. That is, a decorated permutation is a permutation with fixed points colored in two colors.
\end{definition}

It is easy to see the bijection between necklaces and decorated permutations. To go from a Grassmann necklace $\I$ to a decorated permutation $\pi^{:}=(\pi,col)$

\begin{itemize}
 \item if $I_{i+1} = (I_i \backslash \{i\}) \cup \{j\}$, $j \not = i$, then $\pi(i)=j$
 \item if $I_{i+1} = I_i$ and $i \not \in I_i$ then $\pi(i)=i, col(i)=1$
 \item if $I_{i+1} = I_i$ and $i \in I_i$ then $pi(i)=i, col(i)=-1$
\end{itemize}

To go from a decorated permutation $\pi^{:}=(\pi,col)$ to a Grassmann necklace $\I$
$$I_r = \{ i \in [n] | i<_r \pi^{-1}(i) \textbf{ or } (\pi(i)=i \textbf{ and } col(i)=-1) \}.$$
Recall we have defined $<_r$ to be a total order on $[n]$ such that $r<_r r+1 <_r \cdots <_r n <_r 1 <_r \cdots <_r r-1$. This is same as $<_{c^{r-1}}$ where $c=(1,\cdots,n) \in S_n$.

\begin{lemma}[\cite{P}]
\label{lem:P}
For a matroid $\M \subseteq {[n] \choose k}$ of rank $k$ on the set $[n]$, let $\I_{\M} = (I_1,\cdots,I_n)$ be the sequence of subsets such that $I_i$ is the minimal member of $\M$ with respect to $\leq_i$. Then $\I_{\M}$ is a Grassmann necklace.
\end{lemma}

\begin{theorem}[\cite{O}]
\label{thm:O}
$\M$ is a positroid if and only if for some Grassmann necklace $(I_1, \cdots, I_n)$,
$$ \M = \bigcap_{i=1}^{n} SM_{I_i}^{c^{i-1}} $$
In other words, $\M$ is a positroid if and only if the following holds : $H \in \M$ if and only if $H \geq_t I_t$ for any $t\in [n]$.
\end{theorem}

In \cite{P}, Postnikov also showed a bijection between each positroid cells and an combinatorial object called $\LE$-diagrams. So for a $\LE$-diagram $D$, let's denote $\S_{D}$ as the cell in the nonnegative Grassmannian corresponding to $D$. $\LE$-diagrams come from $\LE$-graphs, which are very special cases of plabic-graphs, which we will defined in section 4.

\begin{definition}
A Young diagram of shape $\lambda$ is called a $\LE$-diagram if it satisfies the following property. Each box is either empty or filled with one dot. For any three boxes indexed $(i,j),(i',j),(i,j')$, where $i'<i$ and $j'<j$, if boxes on position $(i',j)$ and $(i,j')$ contain a dot inside, then the box on $(i,j)$ also contains a dot. This property is called the $\LE$-property.
\end{definition}

\begin{figure}
	\centering
		\includegraphics{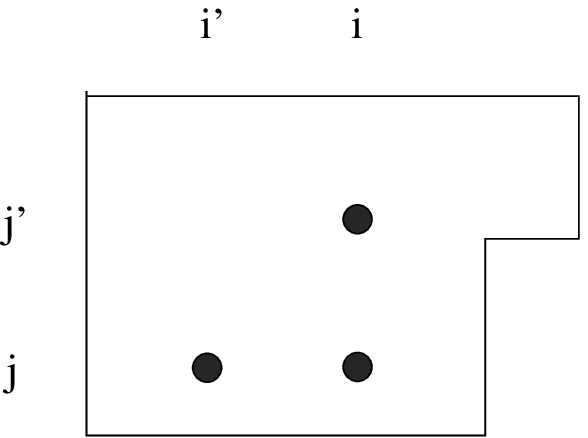}
	
	\caption{$\LE$-property}
	\label{fig:lehcondition}
\end{figure}

\begin{figure}
	\centering
		\includegraphics{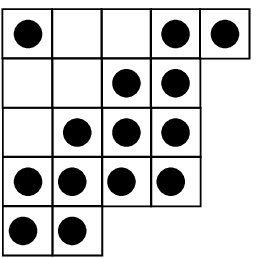}
	
	\caption{Example of a $\LE$-diagram}
	\label{fig:lediagex}
\end{figure}

The boundary of the diagram forms a lattice path from the upper-right corner to lower-left corner. Label the $n$ steps in this path by numbers $1,\cdots,n$ consecutively. Define $I(\lambda)$ as the set of lables of $k$ vertical steps in the path. Put dots on each edge of the boundary path. Connect all dots on same row and connect all dots on the same column. Then we get a $\LE$-graph.

Let's say a dot in the diagram is at $(x,y)$ if it is at $x$-th row and $y$-th column. Denote $I(\lambda)$ simply by $I$. Let $I = \{i_1,\cdots,i_k\}, [n]\setminus I = \{j_1,\cdots,j_{n-k}\}$ such that $i_1 < \cdots < i_k, j_1> \cdots > j_{n-k}$. The Pl\"ucker variable $\Delta_{X,Y}$ is defined as $\Delta_{ (I \setminus \{ i_{x_1},\cdots,i_{x_k} \}) \cup \{j_{y_1},\cdots,j_{y_k} \}}$ where $X = \{x_1,\cdots,x_r\},Y=\{y_1,\cdots,y_r\}, r \leq k$. We define the rank of $\Delta_{X,Y}$ as $|X|=|Y|$. Rank of the cell is defined as the maximal rank among nonzero Pl\"ucker variables. For convenience, we will write $aX$ for $\{a\} \cup X$ where $a < min(X)$. $Xb$ denotes $X \cup \{b\}$ where $b > max(X)$. $aXb$ stands for $X \cup \{a,b\}$ where $a<min(X), b> max(X)$.

\begin{definition}
$\LE$-graph is obtained from a $\LE$-diagram in the following way. Put a dot on center of each edge of the boxes on the southeast boundary of the diagram, and label them $1,2,\cdots$ starting from northeast to southwest. Call these the boundary vertices. Now for each dot inside the $\LE$-diagram, draw a horizontal line to its right, and vertical line to its bottom until it reaches the boundary of the diagram.
\end{definition}

\begin{figure}
	\centering
		\includegraphics{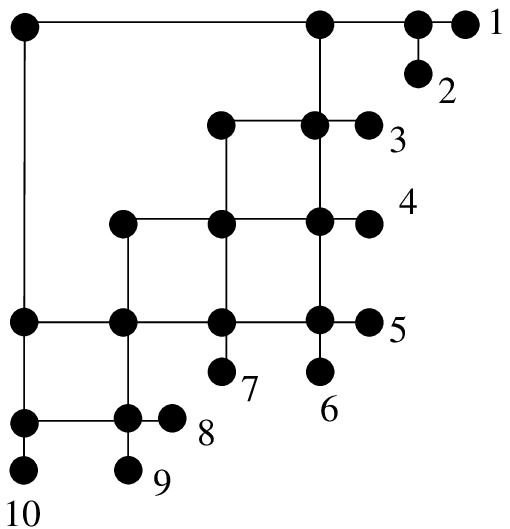}
	
	\caption{Example of a $\LE$-graph}
	\label{fig:legraphex}
\end{figure}

\begin{figure}[ht]
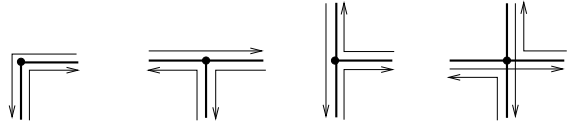
\caption{Rules of the road for $\LE$-graphs}
\label{fig:gamma_rules}
\end{figure}

\begin{theorem}[\cite{P}]
Define a map $\chi$ that sends a $\LE$-diagram to a decorated permutation $\pi^{:}=(\pi,col )$ defined as following. Set $\pi(i)=j$ where we reach $j$ when we start from $i$ and follow the rules of the road in Figure~\ref{fig:gamma_rules}. If $\pi(i)=i$, set $col(i)=-1$ if $i$ is on a horizontal edge, $col(i)=1$ if otherwise. Then $\chi$ is a bijection between $\LE$-diagrams having lower boundary $I$ and decorated permutations having $I_1 = I$, where $\I = (I_1,\cdots,I_n)$ is the Grassmann necklace of the decorated permutation.
\end{theorem}

Lattice path matroids were defined in \cite{BMN}. The name came from the fact that a lattice path matroid represents all lattice paths bounded between two lattice paths. They are special case of positroids.

\begin{definition}
\textit{Lattice path matroids} are defined as the following. Let us be given a base set $[n]$ and $I,J \in {[n] \choose k}$ such that $I \leq J$.
$$ LP_{I,J} = \{H| H \in {[n] \choose k}, I \leq H \leq J \} = SM_I \cap \tilde{SM}_J $$
\end{definition}

Since $I,J$ corresponds to two lattice paths in a $(n-k)$-by-$k$ grid, $LP_{I,J}$ expresses all the lattice paths between them. Lattice path matroid is a positroid.

\begin{theorem}[\cite{O}]
\label{thm:lpm}
Let us be given a base set $[n]$, $I=\{i_1,\cdots,i_k\},J=\{j_1,\cdots,j_k\} \in \ {[n] \choose k}$ such that $i_1<\cdots<i_k,j_1<\cdots<j_k$ and $I \leq J$. Then $LP_{I,J}$ is a positroid and corresponds to the decorated permutation $\pi^{:}=(\pi,col)$ defined as the following.
$$\mbox{$\pi(j_r) = i_r$ for all $r \in [k]$}$$
$$\mbox{$\pi(d_r) = c_r$ for all $r \in [n-k]$}$$
$$ \mbox{If $\pi(t)=t$ then $col(t)$} = \{ \begin{array}{ll}
        -1 & \mbox{if $t \in J$}\\
        1 & \mbox{otherwise} \end{array} $$ 
        
where $[n] \setminus J = \{d_1,\cdots,d_{n-k}\}, [n] \setminus I = \{c_1,\cdots,c_{n-k}\}$ such that $d_1 < \cdots < d_{n-k}, c_1<\cdots<c_{n-k}$.

\end{theorem}

\begin{definition}
A lattice path inside a Young diagram is defined as the following. Start from the southwest corner of the diagram. We follow the edges of the boxes inside the diagram, going up or right all the time till we reach the northeast corner of the diagram.
\end{definition}

\begin{proposition}
\label{prop:lpdiag}
Cells corresponding to lattice path matroids have $\LE$-diagram that looks like the following. Draw a lattice path inside the Young diagram. Then put dots inside the boxes that are to the bottom or to the right of the lattice path. See Figure~\ref{fig:lpex}.

\begin{figure}
	\centering
		\includegraphics{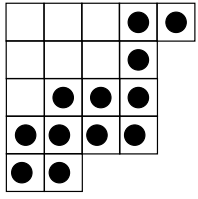}
		\caption{Example of a $\LE$-diagram described in Lemma~\ref{prop:lpdiag}.}
	\label{fig:lpex}
\end{figure}

\end{proposition}
\begin{proof}
Fix any Young diagram. Denote its shape by $\lambda$. Fix $I$ as $I(\lambda)$.
We will prove that decorated permutations coming from cells described in the statement has the form described in Theorem~\ref{thm:lpm}. This with the fact that number of $LP_{I,J}$'s equal the number of lattice paths inside the diagram proves the lemma. 

In the $\LE$-graph, for $i \in I$, call the leftmost dot on the same row the anker point of $i$. See Figure~\ref{fig:lpgrpex}. To determine $\pi^{-1}(i)$, we move south, then east, then south, then east and repeat until we reach $j \in [n]$. This $j$ being $\pi^{-1}(i)$ follows by the rules of the road, since we have walked it in an opposite way. So $\pi^{-1}(i)$ is completely determined by the anker point of $i$. Now for $i'>i$, anker point of $i'$ would place south-west to $i$. (Either strictly south or strictly west or both.) This implies $\pi^{-1}(i') > \pi^{-1}(i)$. Similar argument works for $i \not \in [n]$. So by argument above, the lemma follows.

\begin{figure}
	\centering
		\includegraphics{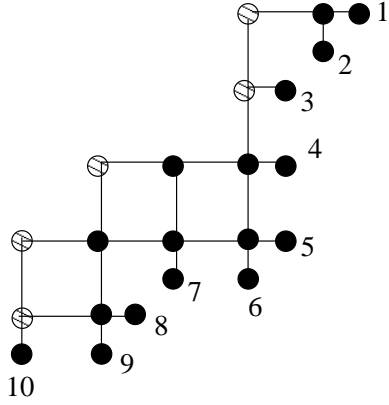}
	
	\caption{Anker points of $i \in I_1$ are highlighted.}
	\label{fig:lpgrpex}
\end{figure}

\end{proof}

\section{Cluster algebra and quivers}
In this section we review the definition of cluster algebras. We first follow the original definition(\cite{FZ}), then express them again in terms of quivers.

A cluster algebra $\A$ is a commutative ring contained in an ambient field $\F$ isomorphic to the field of rational functions in $m$ variables over $\Q$. $\A$ is generated inside $\F$ by a set of generators. These generators are obtained from an initial seed via an iterative process of seed mutations which follows a set of canonical rules.

A seed in $\F$ is a pair $(x,B)$, where
\begin{itemize}
\item $x=\{x_1,\cdots,x_m\}$ is a set of $m$ algebraically independent generators of $\F$, which is split into a disjoint union of an $n$-element cluster $x_c = \{x_1,\cdots,x_n\}$ and an $(m-n)$-element set of frozen variables $x_f = \{x_{n+1},\cdots,x_m\}$
\item $B=(b_{ij})$ is an $m$-by-$n$ integer matrix of rank $n$ whose principal part $B_p = (b_{ij})_{i,j \in [n]}$ is skew-symmetrizable, i.e., there exists a diagonal matrix $D$ with positive diagonal entries such that $DBD^{-1}$ is skew-symmetric.
\end{itemize}
The matrix $B$ is called the exchange matrix of a seed. A seed mutation $\mu_k$ in direction $k \in \{1,\cdots,n\}$ transforms a seed $(x,B)$ into another seed $(x',B')$ defined as follows:
\begin{itemize}
\item $x' = x \setminus \{x_k\} \cup \{x_k'\}$ where $x_k'$ is found from exchange relation $$x_k x_k' = \prod_{b_{ik} >0} {x_i}^{b_{ik}} + \prod_{b_{ik}<0} {x_i}^{-b_{ik}} $$
\item $B' = \mu_k(B)$, i.e., $B$ undergoes a matrix mutation.
\end{itemize}

Now we can express the seeds and seed mutations using quivers as in \cite{K}.

Quiver is an oriented graph. Loop in a quiver is an arrow whose source coincides with its target. 2-cycle is a pair of distinct arrows $\alpha,\beta$ such that the source of $\alpha$ equals target of $\beta$ and vice versa. In this paper, we will restrict ourselves to quivers with no loops and no 2-cycles.

\begin{definition}
A seed is a pair $(x,R)$, where
\begin{itemize}
\item $R$ is a finite quiver without loops or 2-cycles with vertex set $\{1,\cdots,m\}$
\item $x$ is defined as above.
\end{itemize}
Mutation $\mu_k(x,R), 1 \leq k \leq n$ of $(x,R)$ at $k$ is the seed $(x',R')$ where
\begin{itemize}
\item $x' = x \setminus \{x_k\} \cup \{x_k'\}$ where $x_k'$ is found from exchange relation
$$x_k ' x_k = \prod_{i \rightarrow k} x_i + \prod_{k \rightarrow j} x_j $$
\item $R'$ is obtained from $R$ as follows:
\begin{enumerate}
\item Reverse all arrows incident with $k$
\item For all vertices $i \not = j$ distinct from $k$, modify the number of arrows between $i$ and $j$ as following:
\begin{itemize}
\item If there were $s$ edges $i \rightarrow k$, $t$ edges $k \rightarrow j$ and $r$ edges $i \rightarrow j$, then in $R'$, let there be $r+st$ edges $i \rightarrow j$.
\item If there were $s$ edges $k \rightarrow i$, $t$ edges $j \rightarrow k$ and $sr$ edges $i \rightarrow j$, then in $R'$, let there be $r-st$ edges $i \rightarrow j$.
\end{itemize}

\end{enumerate}

\end{itemize}
\end{definition}

Let $Q$ be a finite quiver with vertex set $\{1,\cdots,m\}$. Consider the initial seed $(x=\{x_1,\cdots,x_m\},Q)$.
\begin{itemize}
\item Clusters are defined to be set of $x'$'s appearing in seeds $(x',Q')$ obtained from $(x,Q)$ by mutations described above.
\item Cluster variables are defined to be union of elements of all clusters.
\end{itemize}

\begin{theorem}[\cite{FZ4}]
Any cluster variable is expressed in terms of the variables $x_1,\cdots,x_m$ of any given seed as a Laurent polynomial with integer coefficients.
\end{theorem}

The Laurent phenomenon is itself a very interesting subject, and there has been many interesting research about it. For more, see \cite{FZ3} and \cite{HS}. The following conjecture is called the nonnegativity conjecture for cluster algebras. This is one of the motivations for studying TP-bases.

\begin{conjecture}[\cite{FZ4}]
Every coefficient in these Laurent polynomials is nonnegative.
\end{conjecture}

The above conjecture has been proven for all cluster algebras arising from acyclic quivers(\cite{CR}). For a lot of the cells in the nonnegative Grassmannian, its quiver coming from canonical plabic graph has many cycles. So it is not known whether we can apply the results of \cite{CR} to cluster algebras arising from cells of $Gr_{kn}^{tnn}$. 

A cluster algebra is of finite type if the number of cluster variables are finite (\cite{FZ2}). For cluster algebra of finite type, one of its seeds has its quiver corresponding to a Dynkin quiver. That is, if we delete dots corresponding to the frozen variables and its adjacent edges, and unorient all remaining edges, it becomes a simply laced Dynkin diagram. Aside from the above conjecture being proven, there has been many interesting research on cluster algebra of finite type. In (\cite{CP}), Caroll and Price gave a combinatorial interpretation of the cluster variables in type $A_n$ using paths and perfect matchings. In (\cite{M}), Musiker gave a combinatorial interpretation of cluster variables of type $A_n,B_n,C_n,D_n$ in terms of perfect matchings of graphs. 

\begin{remark}
Cluster algebras of type A,D,E can be constructed from positroid cells. We will show this in the next section, in Proposition~\ref{prop:fintype}.
\end{remark}

In the next section, we will define plabic graphs that corresponds to each cells in the non-negative Grassmannian. The dual graph of plabic graphs will be quivers, so we will be able to define a cluster algebra from those plabic graphs.

\section{Plabic graphs}

\begin{definition}
A planar bicolored graph, or simply a plabic graph is a planar undirected graph $G$ drawn inside a disk. The vertices on the boundary are called boundary vertices. All vertices in the graph are colored either white or black.
\end{definition}

\begin{definition}
For a plabic graph $G$, a trip(one-way trip in notion of \cite{P}) is a directed path $T$ in $G$ such that $T$ joins two boundary vertices and satisfies the following rules of the road. Turn right at a black vertex, and turn left at a white vertex. Trip permutation $\pi_G \in S_n$ is defined such that $\pi_G(i)=j$ whenever trip that starts at the boundary vertex labeled $i$ ends at boundary vertex $j$. Decorated trip permutation $\pi_G^{:}$ is defined similarly.
\end{definition}

We now show how to convert $\LE$-diagrams defined in the previous section to plabic graphs. Then the rules of the road we defined for $\LE$-graphs would coincide with the rules of the road for plabic graphs.

\begin{definition}
First transform the $\LE$-diagram into a $\LE$-graph. Then draw a big circle, so that the boundary vertices of the $\LE$-graph all lie on the circle. Then transform each dots according to the rule in Figure~\ref{fig:gamma_plabic}. We will call this the canonical plabic graph of the cell (The cell in the nonnegative Grassmannian that corresponds to the $\LE$-diagram).

\begin{figure}[ht]
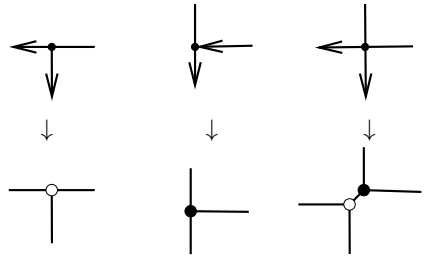
\caption{Transforming $\LE$-graphs into plabic graphs}
\label{fig:gamma_plabic}
\end{figure}

\end{definition}

We have following 3 moves on plabic graphs. In this paper, we will be looking at plabic graphs obtained from a canonical plabic graph by the 3 moves below. So we will be restricting ourselves to reduced plabic graphs in sense of \cite{P}.

\medskip

\begin{figure}[ht]
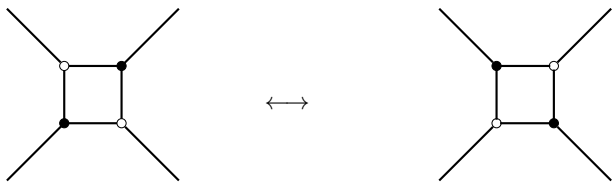
\caption{(M1) Square move}
\label{fig:square_move} 
\end{figure}

(M1) Pick a square with vertices alternating in colors. Then we can switch the colors of all the vertices. See Figure~\ref{fig:square_move}.

\begin{figure}[ht]
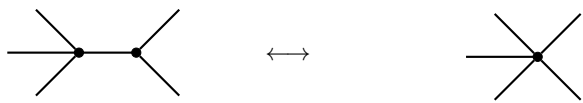
\caption{(M2) Unicolored edge contraction}
\label{fig:edge_contraction}
\end{figure}

(M2) For two adjoint vertices of the same color, we can contract them into one vertice. See Figure~\ref{fig:edge_contraction}.

\begin{figure}[ht]
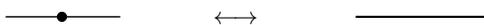
\caption{(M3) Middle vertex insertion/removal}
\label{fig:vertex_removal}
\end{figure}

(M3) We can insert or remove a vertice inside any edge. See Figure~\ref{fig:vertex_removal}.

\begin{theorem}[\cite{P}] Let $G$ and $G'$ be two (reduced) plabic graphs with the same number of boundary vertices. Then the following claims are equivalent:
\begin{itemize}
\item $G$ can be obtained from $G'$ by moves (M1)-(M3).
\item These two graphs have the same decorated trip permutation $\pi_G^{:} = \pi_{G'}^{:}$.
\end{itemize}

\end{theorem}

Now we will show how one can define a cluster algebra for each cell in the nonnegative Grassmannian using the canonical plabic graph. For a plabic graph $G$, look at its dual graph $G'$. Each vertice in $G'$ corresponds to a face of $G$, and two vertices in $G'$ are connected by an edge if and only if the corresponding faces of $G$ share an edge, and that edge connects two vertices of different color. We orient the edges according to the rule below, to obtain a quiver. 

\begin{figure}[ht]
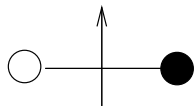
\caption{Giving an orientation to the dual graph}
\label{fig:oridual}
\end{figure}

\begin{figure}
	\centering
		\includegraphics{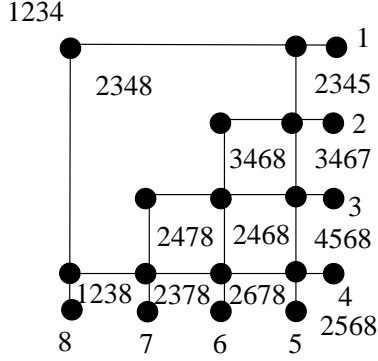}
	
	\caption{Labeling faces of the $\LE$-graph. }
	\label{fig:legraphlabel}
\end{figure}

All we have to do now is to label the faces of the canonical plabic graph with certain variables to define the initial seed. For each trip inside the plabic graph that starts at $i$ and end at $j$, write $j$ inside each face that lies left to that trip (This is well defined since the entire graph is bounded by a circle, and $i,j$ lie on the circle). For a face written $i_1, \cdots, i_k$, label the face with the Pl\"ucker variable $\Delta_{i_1 \cdots i_k}$. For $\LE$-graphs, the labeling is defined similarly since the rules of the road coincides with that of plabic graphs. See Figure~\ref{fig:legraphlabel}. Set the variables corresponding to faces on the boundary to be frozen variables. We also call these frozen variables as variables of the Grassmann necklace. So we have finished defining a cluster algebra associated to each cell in the nonnegative Grassmannian. Call them the canonical cluster algebra of the cell.

\medskip

There is a simpler way to read off the Pl\"ucker variables of each faces in the canonical plabic graph. Inside the $\LE$-diagram, pick any dot $x$ at position $(a,b)$. Then pick a dot $y$ as the following way:

Look at dots inside region $\{(x,y)|x<a, y<b\}$. Give a partial ordering such that for dots of same row, dot to the right is bigger, and for dots of same column, dot below is bigger. Then there is a unique maximal dot in the region due to the $\LE$-condition. This will be our $y$.

\begin{figure}
	\centering
		\includegraphics{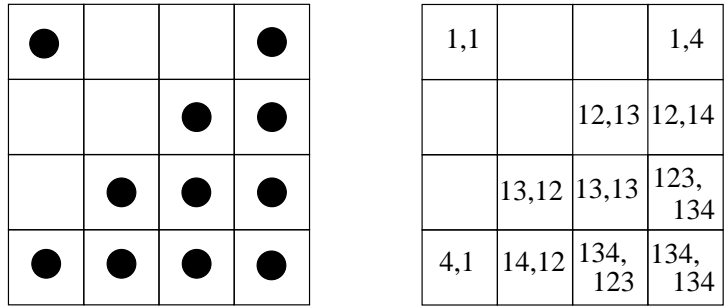}
	\label{fig:lediaglabel}
	\caption{Labeling faces of the $\LE$-diagram. $\{I,J\}$ denotes $\Delta_{I,J}$. }
\end{figure}

If we found a dot $y$, then we will say that $y$ covers $x$. If there is no such dot, we will say that $x$ is uncovered. For a dot $x_1$, pick a sequence of dots $x_2,\cdots,x_t$ such that $x_{i+1}$ covers $x_i$ for $ 1 \leq i \leq t-1$ and $x_t$ is uncovered. Then we label the face that corresponds to the dot $x_1$ by $\Delta_{\{a_1 \cdots a_t\}, \{b_1 \cdots b_t\}}$ where $(a_i,b_i)$ is the position of $x_i$ for $1 \leq i \leq t$. See Figure~\ref{fig:lediaglabel}.

\begin{proposition}[PW]
\label{lem:PW}
The labeling for $\LE$-diagram above is compatible with the labeling of canonical plabic graph.
\end{proposition}

Now as we have promised in the previous section, let's prove that we can construct finite-type cluster algebras of type $A_n,D_n,E_n$. 

\begin{proposition}
\label{prop:fintype}
There are positroid cells such that the associated cluster algebra is of type $A_n,D_n,E_n$. 
\end{proposition}

\begin{proof}
Let us denote a $\LE$-diagram by $D$, its plabic graph by $G_D$ and quiver obtained from dual graph of $G_D$ by $Q_D$. Two vertices in $Q_D$ are connected if and only if the corresponding dots $x,y$ in $D$ satisfy one of the following.

\begin{enumerate}
\item $x,y$ are adjacent to each other and is on the same row.
\item $x,y$ are adjacent to each other and is on the same column.
\item $x$ covers $y$ or $y$ covers $x$.
\end{enumerate}

So we need to find $\LE$-diagrams such that when we delete all the dots corresponding to Grassmann necklace and connect the dots that satisfy one of the three conditions above, we get a simply laced Dynkin diagram. The rest is trivial. We show how to construct $A_7,D_7$ in Figure~\ref{fig:fintype}. The rest of $A_n,D_n,E_n$ are constructed similarly.

\begin{figure}
	\centering
		\includegraphics{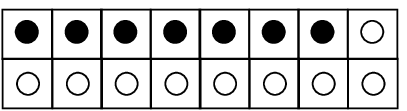}
		\includegraphics{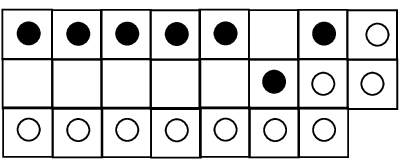}
	\label{fig:fintype}
	\caption{$\LE$-diagrams corresponding to cluster algebra of type $A_7,D_7$. White dots correspond to dots of the Grassmann necklace. }
\end{figure}

\end{proof}

\section{$\LE$-networks}

In this section we explain how $\LE$-graphs are used to parametrize the positroid cells. For a $\LE$-graph, direct each horizontal edge to the left, and vertical edge to the right. Call vertices labeled with numbers as the boundary vertices, and all other vertices as internal vertices. Put a weight on each edge, but let the value be $1$ for all vertical edges. We call this a $\LE$-network.

\begin{figure}
\includegraphics[height=1.2in]{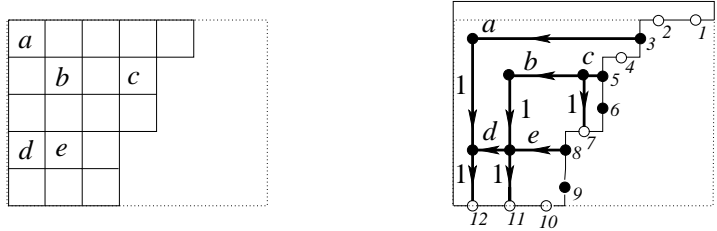}
\caption{A $\LE$-tableau and the corresponding $\LE$-network}
\label{fig:le_tableau_gamma}
\end{figure}


For a $\LE$-network with $n-k$ sinks and $k$ sources, define the boundary measurement matrix $A \in Mat_{kn}$ as the following. Denote $I$ as the set of indices of sources, $I^c$ as the set of indices of sinks.
\begin{itemize}
\item Submatrix $A_I$ is the identity matrix $Id_k$.
\item Set $M_{ij}$ as the following.
$$M_{ij} := \sum_{P:b_i \rightarrow b_j} \prod_{e \in P} x_e$$
Remaining entries of $A$ are $a_{rj} = (-1)^{s} M_{i_r,j}$ for $r \in [k]$ and $j \in I^c$ where $s$ is the number of elements of $I$ strictly between $i_r$ and $j$.
\end{itemize}
Then we have the following result.

\begin{theorem}[\cite{P}]

\begin{enumerate}
\item $\Delta_{H,J}(A)/\Delta_{\phi,\phi}(A) = \sum_{P: H \rightarrow J} \prod_{e \in P} x_e$ where $P: I \rightarrow J$ indicates non-vertice-crossing paths from source $H$ to sinks $J$. 
\item Matrix $A$ is parametrized by $x_{e_1},\cdots,x_{e_d}$ where $e_1,\cdots,e_d$ stands for horizontal edges of the $\LE$-network. $d$ will also be same as the number of dots. So define a map $\phi$, that sends $(x_{e_1},\cdots,x_{e_d})$ to this $A$. Then if this $\LE$-network came from positroid $\M$, $\phi( \R_{>0}^{d}) = S_M^{tnn}$.
\end{enumerate}

\end{theorem}

If we represent the elements of $Gr_{kn}$ as square matrices of size $(n-k)$ by ignoring columns $I$, the first statement of the theorem can be thought of as the Lindstrom lemma (\cite{LI}).

The above theorem tells us that we can express each Pl\"ucker variables of the cells as polynomials in horizontal-edge-weights and $\Delta_{\phi,\phi}$.

\begin{theorem}[\cite{P}]
Values of Pl\"ucker variables coming from labeling of $\LE$-diagram (including $\Delta_{\phi,\phi})$ completely determines values of all other Pl\"ucker variables. Which also means that it completely determines values of all horizontal-edge-weights.
\end{theorem}

We will say that a variable $\Delta_{H,J}$ is a unique-path-variable in a $\LE$-network (or a cell corresponding to that $\LE$-network) if there is only one family of non-vertice-crossing paths $H \rightarrow J$.

The following lemma is why we want to transform a TP-bases to a TP-bases consisting of unique-path-variables of the cell.

\begin{lemma}
\label{lem:uniqpath}
Assume we have a TP-bases $\S$ such that all of its variables corresponds to unique paths in the $\LE$-network. Then all other Pl\"ucker variables can be expressed in a Laurent polynomial with positive coefficients in variables of $\S$.
\end{lemma}
\begin{proof}
All weights of the edges can be expressed as Laurent monomials in variables of $\S$. 
\end{proof}

\section{TP-bases and TP-diagrams}

We will define TP-bases and TP-diagrams in this section. From now on, we will use the following terms interchangably.

\begin{itemize}
\item Cell in the nonnegative Grassmannian
\item $\LE$-diagram of the cell
\item $\LE$-graph of the cell
\item $\LE$-network of the cell
\end{itemize}

First let's define the TP-bases. We define TP-bases for each cell in the nonnegative Grassmannian.

\begin{definition}
Let's define TP-bases as the set of algebraically independent Pl\"ucker variables $\B$ such that all other nonzero variables in the cell can be expressed as subtraction-free rational function in variables of $\B$.
\end{definition}

\begin{definition}
We say that we obtained set of Pl\"ucker variables $\S'$ from set of Pl\"ucker variables $\S$, by Pl\"ucker mutation if one of the following holds. For a 3-term Pl\"ucker relation
$$\Delta_{H_1}\Delta_{H_2} = \Delta_{H_3}\Delta_{H_4} + \Delta_{H_5}\Delta_{H_6}$$ 
\begin{itemize}
\item $\Delta_{H_i} \in \S$ for $2 \leq i \leq 6$. $\Delta_{H_1} \not \in \S$. $\S' = (\S \setminus \{\Delta_{H_2}\}) \cup \{\Delta_{H_1}\}$.
\item $\Delta_{H_3}$ or $\Delta_{H_4}$ is zero. $\Delta_{H_5},\Delta_{H_6},\Delta_{H_2} \in \S$. $\Delta_{H_1} \not \in \S$. $\S' = (\S \setminus \{\Delta_{H_2}\}) \cup \{\Delta_{H_1}\}$.

\item $\Delta_{H_3}$ or $\Delta_{H_4}$ is zero. $\Delta_{H_1},\Delta_{H_2},\Delta_{H_6} \in \S$. $\Delta_{H_5} \not \in \S$. $\S' = (\S \setminus \{\Delta_{H_6}\}) \cup \{\Delta_{H_5}\}$.
\end{itemize}
\end{definition}

\begin{remark}
\label{rem:mutationbases}
Assume $\S$ was the TP-bases of some positroid cell. Let us obtain $\S'$ by a Pl\"ucker mutation. Then $\S'$ is also a TP-bases.
\end{remark}

The set of Pl\"ucker variables obtained from labeling of $\LE$-diagram forms a TP-bases. We want to show that set of Pl\"ucker variables obtained from $\LE$-diagram can be transformed using Pl\"ucker mutations, to a set of Pl\"ucker variables that corresponds to unique paths in the $\LE$-graph. Then by Lemma~\ref{lem:uniqpath}, Conjecture~\ref{con:main} follows. Now let's define the TP-diagram.

\begin{definition}
A pre-TP-diagram is a Young diagram with set of Pl\"ucker variables. For each variable $\Delta_{Xx,Yy}$, we put a dot in the diagram at box $(x,y)$, which stands for $x$-th row and $y$-th column. We will always say that $\Delta_{\phi,\phi}$ is inside the diagram, although it doesn't correspond to any dots.
\end{definition}

Since each dot stands for a variable and each variable corresponds to a dot, we will use the term dots and Pl\"ucker variables interchangeably. 

\begin{definition}
We define the rank of a dot to be $|X|$ where the variable corresponding to the dot is $\Delta_{X,Y}$. Pick two dots $\alpha,\beta$ at $(x,y),(x',y')$. We will say that
\begin{itemize}
\item $\beta$ is right of $\alpha$ if $x=x', y' > y$.
\item $\beta$ is left of $\alpha$ if $x=x', y' < y$.
\item $\beta$ is below $\alpha$ if $x' <x, y'=y$.
\item $\beta$ is above $\alpha$ if $x' >x, y'=y$.
\item $\beta$ is adjacent to $\alpha$ if $\beta$ is right to or left to or below or above $\alpha$ and there are no dots between them.
\item $\alpha$ dominates $\beta$ if $\alpha$ is $\Delta_{\{x_1,\cdots,x_k\},\{y_1,\cdots,y_k\}}$ and $\beta$ is $\Delta_{\{x_1,\cdots,x_{k'}\},\{y_1,\cdots,y_{k'}\}}$, $k<k'$. 
\item Region left(right) to $\alpha$ is defined as set of all dots $(x'',y'')$ such that $y'' < y$($y'' > y$). Region below(above) $\alpha$ is defined as set of all dots $(x'',y'')$ such that $x'' > x$($x'' < x$).
\end{itemize}

We define the $k$-rank-line (Will be abbreviated as $k$-rkline) to be the following. Pick all dots of rank $k$. Connect all such adjacent dots with lines. Although this may not be connected and may not look like a line we will call it the $k$-rkline of the pre-TP-diagram. Rank of the diagram is defined to be maximal value among the rank of the dots inside.

For a dot $\alpha$, we can look at dots inside $k$-rkline dominated by $\alpha$. Call that part of the $k$-rkline as $k/\alpha$-rkline. For a dot, if it has a adjacent dots to the below and right on the same rkline, call it a northwest(NW)-corner. Define the southeast(SE)-corner similarly. See Figure~\ref{fig:corners}. Even if $\beta$ is a SW(NE)-end of a $k/\alpha$-rkline, we still call $\beta$ a NW-corner of $k/\alpha$-rkline.

\begin{figure}[ht]
\includegraphics[height=1.2in]{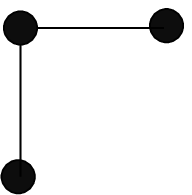}
\includegraphics[height=1.2in]{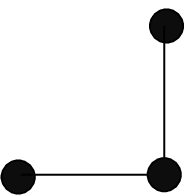}
\caption{Northwest corner and Southeast corner of a rank-line}
\label{fig:corners}
\end{figure}

\end{definition}

\begin{definition}
We call a pre-TP-diagram smooth if it satisfies the following two properties.
\begin{itemize}
\item For any two dots $\gamma_1,\gamma_2$ on $k$-rkline, the following holds. Denote position of $\gamma_i$ by $(x_i,y_i)$. If $x_1 \leq x_2$ then $y_1 \geq y_2$. So we can choose Southwest(SW)-end and Northeast(NE)-end of the $k$-rkline.
\item Adjacent variables on the same rkline differ by one element as the following. If $\Delta_{X_1,Y_1}, \Delta_{X_2,Y_2}$ are on the same rkline and are adjacent, either $X_1=X_2, |Y_1 \setminus Y_2| = 1$ or $|X_1 \setminus X_2|=1, Y_1=Y_2$.

The first condition implies that we can define NW-corners and SE-corners of a rkline.

\end{itemize}

We call a pre-TP-diagram supersmooth if it satisfies the following property. 
\begin{itemize}
\item Pick any two adjacent variables $\Delta_{X,Y},\Delta_{X',Y'}$, $|X| = |X'|-1$. Then $X \subset X', Y \subset Y'$.
\end{itemize}

\end{definition}

\begin{definition}
We call a dot stable if the following properties hold. Let the dot be $\Delta_{\{x_1,\cdots,x_k\},\{y_1,\cdots,y_k\}}$. Look at dots inside region $\{(x,y)|x<x_k, y<y_k\}$. Give a partial ordering such that for dots of same row, dot to the right is bigger, and for dots of same column, dot below is bigger. Then there is a unique maximal dot in the region. That dot is $\Delta_{\{x_1,\cdots,x_{k-1}\},\{y_1,\cdots,y_{k-1}\}}$.

We call a diagram stable if all dots are stable. Diagram being smooth stable that for dot $\Delta_{\{x_1,\cdots,x_k\},\{y_1,\cdots,y_k\}}$, we have $\Delta_{\{x_1,\cdots,x_i\},\{y_1,\cdots,y_i\}}$ inside the diagram for all $1 \leq i \leq k$.
\end{definition}

It is straight from the definitions that $\LE$-diagram is a supersmooth, stable TP-diagram. Now let us define two operations on a pre-TP-diagram, fold and push. Let's first define the fold.

\medskip
Let us pick three adjacent dots of the same rank, of form $\Delta_{X,Y},\Delta_{X,Y'},\Delta_{X,Y'}$, $|X \setminus X'|=|Y \setminus Y'|=1$. In other words, we pick 3 dots of a NW or SE-corner of a rkline. We call the following Pl\"ucker mutation the fold. 
$$\Delta_{X,Y}\Delta_{X',Y'} = \Delta_{X',Y}\Delta_{X,Y'} + \Delta_{X \cap X',Y \cap Y'}\Delta_{X \cup X',Y \cup Y'}$$

We will say that we fold at $\Delta_{X,Y}$. Or if there is no confusion, fold at $(x,y)$ where $x = max(X), y= max(Y)$. Fold denoted above moves a dot at $(x,y)$ to $(x',y')$ where $x' = max(X'), y' = max(Y')$.

\medskip

Now let's define the push. Let us pick any dot of rank $k$ of form $\Delta_{aX,Yb}$ such that the adjacent dot to above is $\Delta_{aX',Yb}$,$|X \setminus X'|=1$ and adjacent dot to the left is $\Delta_{X,Y}$. We call the following Pl\"ucker mutation the push/pull. A push when we replace $\Delta_{aX,Yb}$ with $\Delta_{X',Y}$. A pull when we replace $\Delta_{X',Y}$ with $\Delta_{aX,Yb}$.

$$\Delta_{aX,Yb}\Delta_{X',Y} = \Delta_{X,Y}\Delta_{aX',Yb} + \Delta_{X \cup X', Yb}\Delta_{\{a\} \cup (X \cap X'), Y}$$

Similarly, let us pick any dot of rank $k$ of form $\Delta_{Xa,bY}$ such that the adjacent dot to left is $\Delta_{Xa,bY'}$,$|Y \setminus Y'|=1$ and adjacent dot to the above is $\Delta_{X,Y}$. We call the following Pl\"ucker mutation the push/pull. A push when we replace $\Delta_{Xa,bY}$ with $\Delta_{X,Y'}$. A pull when we replace $\Delta_{X,Y'}$ with $\Delta_{Xa,bY}$.

$$\Delta_{Xa,bY}\Delta_{X,Y'} = \Delta_{X,Y}\Delta_{Xa,bY'} + \Delta_{Xa,Y \cup Y'}\Delta_{ X, \{b\} \cup (Y \cap Y') }$$

When we need to distinguish between two pushes(pulls) defined above, we will denote the first one by SW-push(pull) and the second one by NE-push(pull). We will give plenty of examples of fold and push in the next section.

\begin{remark}
\label{rem:pathrklim}
For all pre-TP-diagrams obtained from a $\LE$-diagram by pulls/pushes and folds, there cannot exist a dot with rank bigger than the rank of the cell.
\end{remark}
\begin{proof}
Fold does not change the rank of the dot being folded. Push/pull does not increase the maximal rank of the dots involved in push/pull.
\end{proof}

\begin{definition}
We say that SW(NE)-end of $k/\alpha$-rkline and $(k+1)/\alpha$-rkline are connected if one of the following holds.
\begin{enumerate}
\item SW(NE)-end of $k/\alpha$-rkline dominates SW(NE)-end of $(k+1)/\alpha$-rkline
\item SW(NE)-end of $k/\alpha$-rkline is adjacent to SW(NE)-end of $(k+1)/\alpha$-rkline
\end{enumerate}
We say that $k/\alpha$-rkline and $(k+1)/\alpha$-rkline are connected if their SW-ends and NE-ends are connected. We call a cell or a $\LE$-diagram associated to it connected if for any dot $\alpha$ and $1 \leq k \leq m-1$ the following two conditions are satisfied.

\begin{itemize}
\item $k$-rkline is connected. That is, $k$-rkline consists of one connected component.
\item $k/\alpha$-rkline is connected with $(k+1)/\alpha$-rkline. 
\end{itemize}

We call a cell or a $\LE$-diagram associated to it weakly connected if the following two conditions are satisfied.
\begin{itemize}
\item $k$-rkline is connected for each $k$.
\item For any dot $\alpha$, one of the following holds.

\begin{itemize}
\item For odd $k$, NW-ends of $k/\alpha$-rkline and $(k+1)/\alpha$-rkline are connected. For even $k$, SE-ends of $k/\alpha$-rkline and $(k+1)/\alpha$-rkline are connected.
\item For odd $k$, SE-ends of $k/\alpha$-rkline and $(k+1)/\alpha$-rkline are connected. For even $k$, NW-ends of $k/\alpha$-rkline and $(k+1)/\alpha$-rkline are connected.
\end{itemize}

\end{itemize}

\end{definition}

\begin{remark}
SW(NE)-end of $k/\alpha$-rkline and $(k+1)/\alpha$-rkline are connected if and only if their corresponding faces in the plabic graph are adjacent.

\end{remark}

\begin{remark}
The cells in Remark~\ref{rem:dbruhatcell} are weakly connected. So our Theorem~\ref{thm:main} is a generalization of a result in \cite{FZ6}.
\end{remark}

\begin{remark}
The cells corresponding to cluster algebras of finite type constructed in the previous section are all connected cells.

\end{remark}

It is easy to see that $\LE$-diagram of any lattice path matroid cells is connected. So when we look at a set of all connected cells, it contains all lattice path matroid cells. Now here is our main theorem.

\begin{theorem}
\label{thm:main}
For weakly connected positroid cells, the canonical TP-bases can be mutated to a TP-bases consisting of unique-path-variables.
\end{theorem}

To prove this theorem, we will use the following Lemma.

\begin{lemma}
\label{lem:main}
Pick any weakly connected positroid cell. Let $\S$ be the set of Pl\"ucker variables read from the $\LE$-diagram. Pick any NW-corner of the 1-rkline. Let $\S'$ be the set of Pl\"ucker variables read from the $\LE$-diagram by deleting that corner. Then we can use Pl\"ucker mutation to get $\S'$ and one variable corresponding to a unique path of $\S$.
\end{lemma}

\section{Proof of Lemma~\ref{lem:main} I : Examples}

Out goal is to prove Lemma~\ref{lem:main}. We will prove the lemma in next three sections. We will be constructing an algorithm to prove the lemma. We will first implicitly show how the algorithm works with a simple example. Below the description of each step, we will write down the Pl\"ucker mutations we have used.

(T1) Let's look at a cell inside $Gr_{4,8}^{tnn}$. 

\begin{figure}[ht]
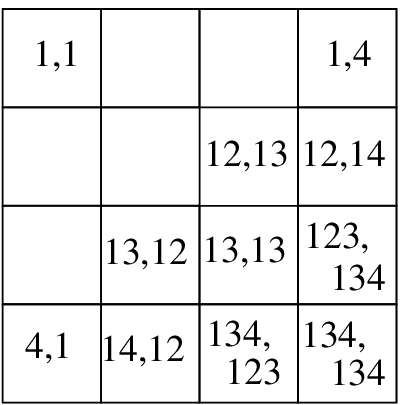
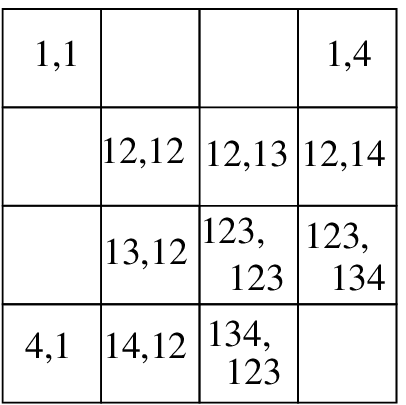
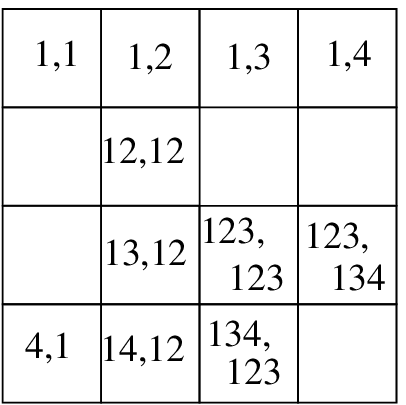
\caption{(T1),(T2),(T3).}
\label{fig:letrasf-1}
\end{figure}

(T2) Apply fold at $\Delta_{134,134}$ and then at $\Delta_{13,13}$. Then all rklines will have no SE-corner. We will call this process the straightening process. We will denote the straightened $\LE$-diagrams as $S\LE$-diagrams.
$$ \Delta_{134,134} \Delta_{123,123} = \Delta_{123,134} \Delta_{134,123} $$
$$ \Delta_{13,13} \Delta_{12,12} = \Delta_{12,13} \Delta_{13,12} + \Delta_{1,1} \Delta_{123,123}$$

(T3) We want to exchange $\Delta_{1,1}$ with $\Delta_{4,2}$. To do so, we need $\Delta_{1,2}$. So exchange $\Delta_{12,14}$ with $\Delta_{1,3}$. Then we can exchange $\Delta_{12,13}$ with $\Delta_{1,2}$.
$$\Delta_{12,14} \Delta_{1,3} = \Delta_{12,13} \Delta_{1,4}$$
$$\Delta_{12,13} \Delta_{1,2} = \Delta_{12,12} \Delta_{1,3}$$

\begin{figure}[ht]
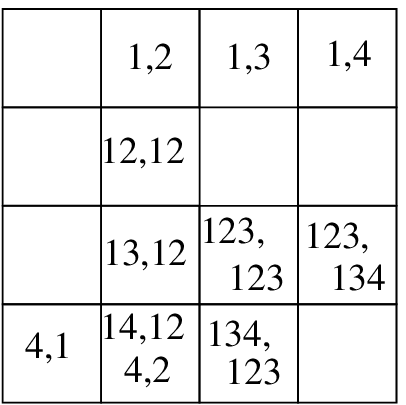
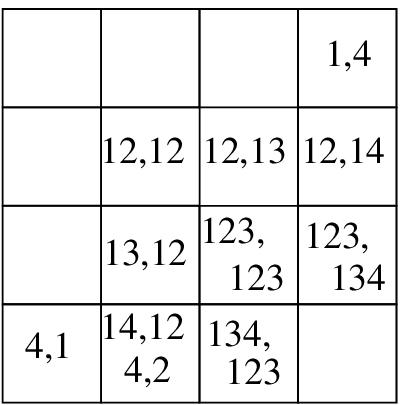
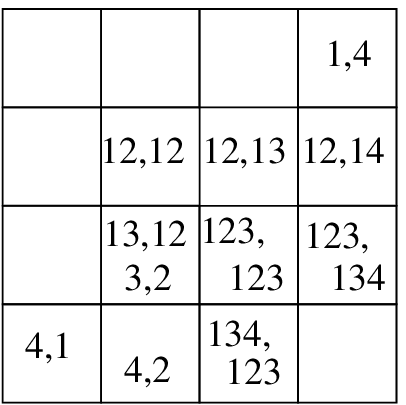
\caption{(T4),(T5),(T6).}
\label{fig:letrasf-2}
\end{figure}

(T4) Now we can exchange $\Delta_{1,1}$ with $\Delta_{4,2}$. We will call this jumping from $1$-rkline to $2$-rkline. There are two ends in the $2$-rkline. Call them SW-end and NE-end. $(4,2)$ is at the SW-end, so call this SW-jump.
$$\Delta_{1,1}\Delta_{4,2} = \Delta_{4,1}\Delta_{1,2} + \Delta_{\phi,\phi}\Delta_{14,12}$$

(T5) Return $\Delta_{1,2}$ to $\Delta_{12,13}$ and $\Delta_{1,3}$ to $\Delta_{12,14}$. 
$$\Delta_{12,14} \Delta_{1,3} = \Delta_{12,13} \Delta_{1,4}$$
$$\Delta_{12,13} \Delta_{1,2} = \Delta_{12,12} \Delta_{1,3}$$

(T6) We start from one end of the 2-rkline. We have two variables inside the box, one rank-1 variable $\Delta_{4,2}$ and one rank-2 variable $\Delta_{14,12}$. We can exchange $\Delta_{14,12}$ with $\Delta_{3,2}$. Notice that $\Delta_{13,12}$ was adjacent to $\Delta_{14,12}$.
$$\Delta_{14,12} \Delta_{3,2} = \Delta_{13,12} \Delta_{4,2}$$

(T7) We repeat the process, replacing $\Delta_{13,12}$ with $\Delta_{2,2}$, $\Delta_{12,12}$ with $\Delta_{2,3}$, $\Delta_{12,13}$ with $\Delta_{2,4}$. So with (T6), we have exchanged $$\Delta_{14,12},\Delta_{13,12},\Delta_{12,12},\Delta_{12,13},\Delta_{12,14},\Delta_{4,2}$$ with $$\Delta_{4,2},\Delta_{3,2},\Delta_{2,2},\Delta_{2,3},\Delta_{2,4},\Delta_{12,14}$$ 
We will call this process sweeping through the 2-rkline.

$$\Delta_{13,12} \Delta_{2,2} = \Delta_{12,12} \Delta_{3,2}$$
$$\Delta_{12,12} \Delta_{2,3} = \Delta_{12,13} \Delta_{2,3}$$
$$\Delta_{12,13} \Delta_{2,4} = \Delta_{12,14} \Delta_{2,4}$$

(T8) We can exchange $\Delta_{12,14}$ with $\Delta_{23,34}$. We will call this jumping from $2$-rkline to $3$-rkline. This is a NE-jump.
$$\Delta_{12,14}\Delta_{23,34} = \Delta_{2,4} \Delta_{123,134}$$

(T9) Then sweep through the 3-rkline. Then only one variable is of form $\Delta_{1X,1Y}$ in the TP-diagram. Denote this diagram as $J\LE$-diagram.
$$\Delta_{123,134} \Delta_{23,23} = \Delta_{123,123} \Delta_{23,34}$$
$$\Delta_{123,123} \Delta_{34,23} = \Delta_{134,123} \Delta_{23,23}$$

\begin{figure}[ht]
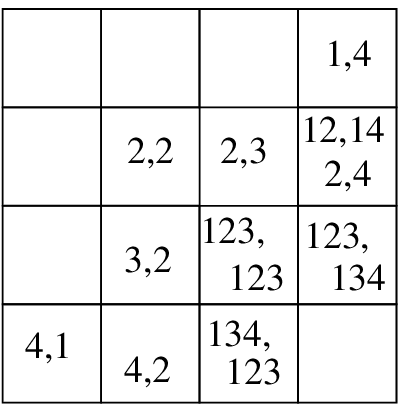
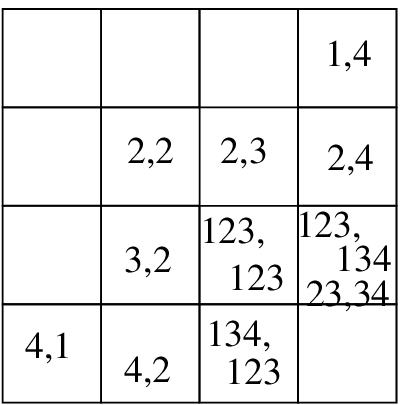
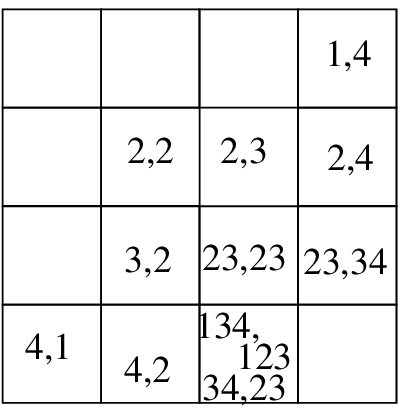
\caption{(T7),(T8),(T9).}
\label{fig:letrasf-3}
\end{figure}

(T10) Fold at $(2,2)$.

(T11) Fold at $\Delta_{23,23}$. Then the resulting TP-diagram looks like a $\LE$-diagram of another cell except for the variable $\Delta_{134,123}$. This variable is a unique-path-variable of the cell. Call this diagram a $G\LE$-diagram.

\begin{figure}[ht]
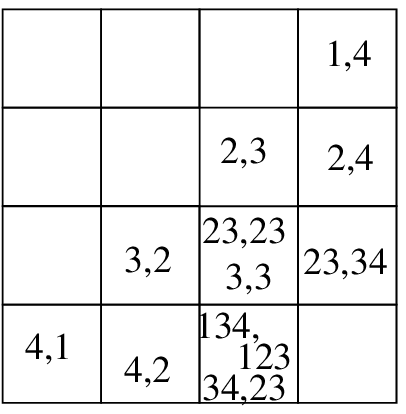
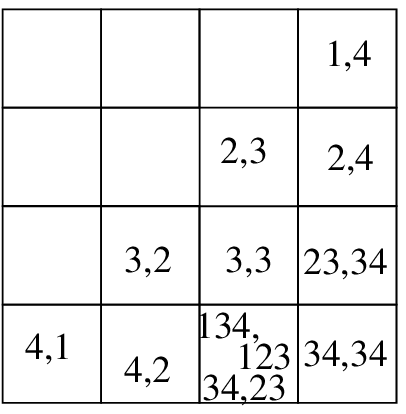
\caption{(T10),(T11).}
\label{fig:letrasf-4}
\end{figure}

We can see that we have obtained a set of variables corresponding to the cell obtained by deleting the NW-corner, and an extra variable that corresponds to a unique path. Now one could think that all we do is straighten the rklines and do the sweeping process. But it is not as simple for other cells, so we have to add another process before the sweeping. Let's look at an example.

(T1) Let's look at a cell contained in inside $Gr_{4,9}^{tnn}$.

(T2) Straighten the rklines, to get a $S\LE$-diagram.

Now assume we exchanged $\Delta_{1,1}$ with $\Delta_{3,2}$ by imitating what we did in (T3),(T4) of the previous example. But then we run into a problem. We can't run the sweeping process since $\Delta_{23,12} \not = 0$, and 
$$\Delta_{13,12} \Delta_{2,2} + \Delta_{23,12}\Delta_{1,2} = \Delta_{12,12} \Delta_{3,2} $$
tells us that it will never be possible to exchange $\Delta_{13,12}$ with $\Delta_{2,2}$ in this cell. So we want $\Delta_{13,12}$ to not be in the 2-rkline for the sweeping process.

\begin{figure}[ht]
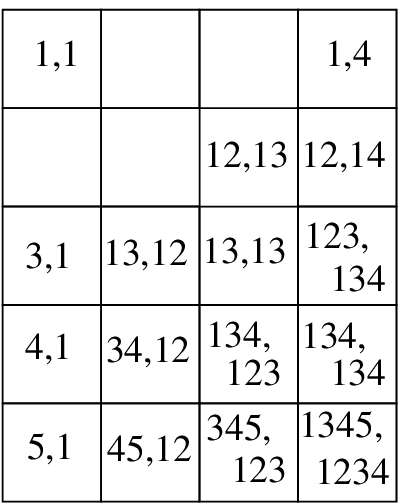
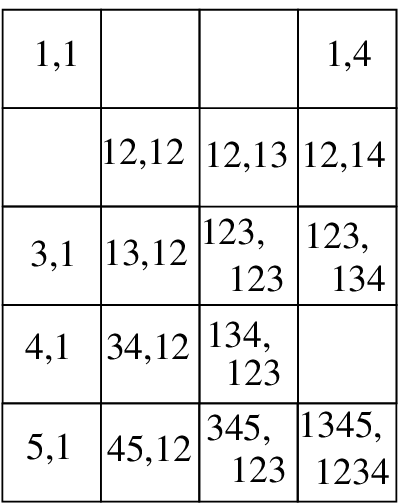
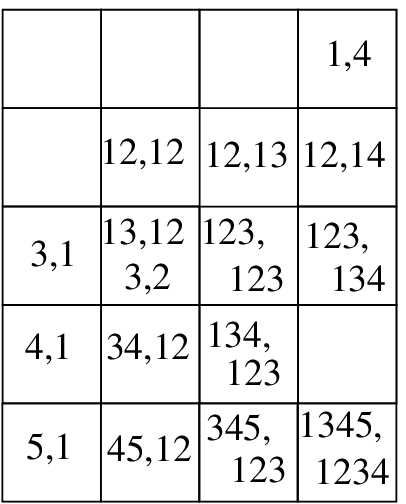
\caption{(T1),(T2). Third figure tells us we can't just go into the sweeping procedure.}
\label{fig:letrasf2-1}
\end{figure}

(T3) Push at $\Delta_{134,123}$ and push at $\Delta_{13,12}$. This means the following.
\begin{enumerate}
\item Exchange $\Delta_{134,123}$ with $\Delta_{23,12}$.
\item Exchange $\Delta_{13,12}$ with $\Delta_{2,1}$. 
\end{enumerate}
We call (T3),(T4) the refining process, which allows us to do sweeps on each rklines in both directions (SW to NE or NE to SW).
$$\Delta_{134,123} \Delta_{23,12} = \Delta_{123,123}\Delta_{34,12}$$
$$\Delta_{13,12} \Delta_{2,1} = \Delta_{12,12}\Delta_{3,1} + \Delta_{23,12}\Delta_{1,1}$$

(T4) Pull back $\Delta_{134,123}$.

$$\Delta_{134,123} \Delta_{23,12} = \Delta_{123,123}\Delta_{34,12}$$

(T5) Push $\Delta_{12,13}$ and $\Delta_{12,14}$, since we want to jump from $\Delta_{1,1}$ to $\Delta_{12,12}$. That is, exchange $\Delta_{1,1}$ with $\Delta_{2,2}$.

$$\Delta_{12,14} \Delta_{1,3} = \Delta_{12,13} \Delta_{1,4}$$
$$\Delta_{12,13} \Delta_{1,2} = \Delta_{12,12} \Delta_{1,3}$$



\begin{figure}[ht]
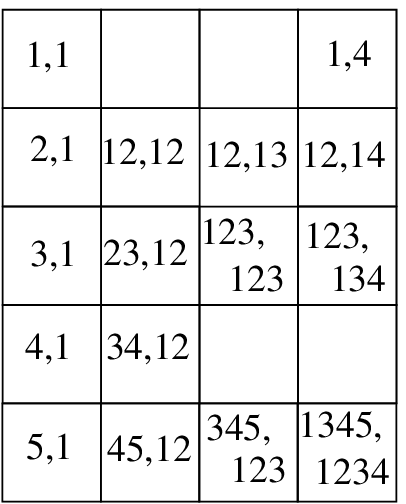
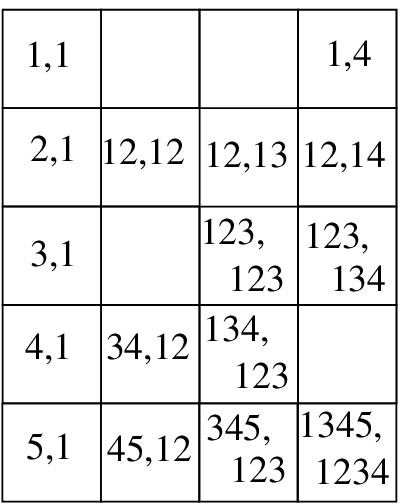
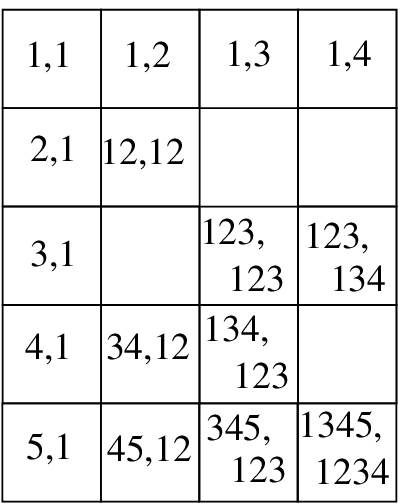
\caption{(T3),(T4), (T5) }
\label{fig:letrasf2-2}
\end{figure}

\begin{figure}[ht]
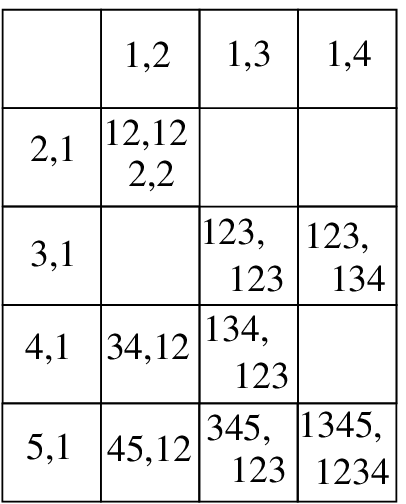
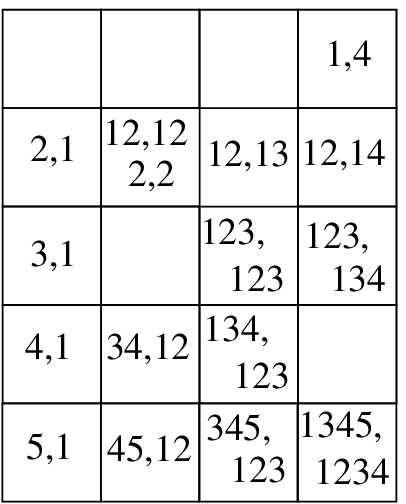
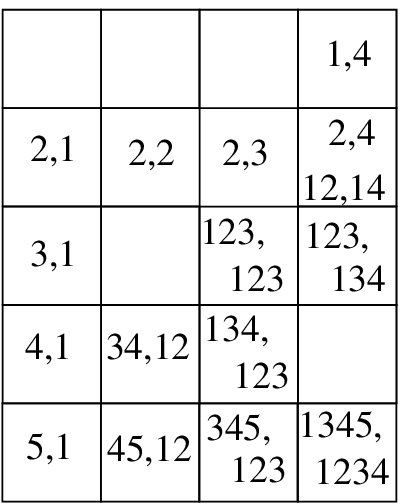
\caption{(T6),(T7),(T8) }
\label{fig:letrasf2-3}
\end{figure}

(T6) Exchange $\Delta_{1,1}$ with $\Delta_{2,2}$.

$$\Delta_{1,1}\Delta_{2,2} = \Delta_{2,1}\Delta_{1,2} + \Delta_{\phi,\phi}\Delta_{12,12}$$

(T7) Pull back $\Delta_{12,13},\Delta_{12,14}$.
$$\Delta_{12,14} \Delta_{1,3} = \Delta_{12,13} \Delta_{1,4}$$
$$\Delta_{12,13} \Delta_{1,2} = \Delta_{12,12} \Delta_{1,3}$$

(T8) Sweep through the $2$-rkline.



\begin{figure}[ht]
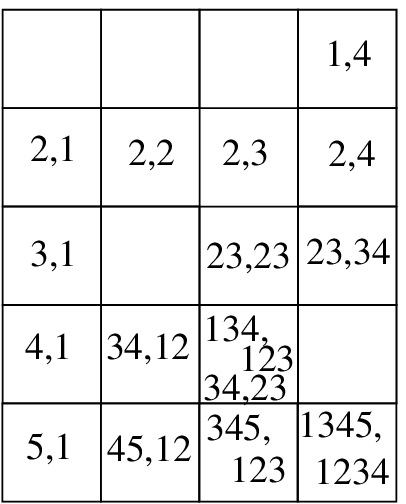
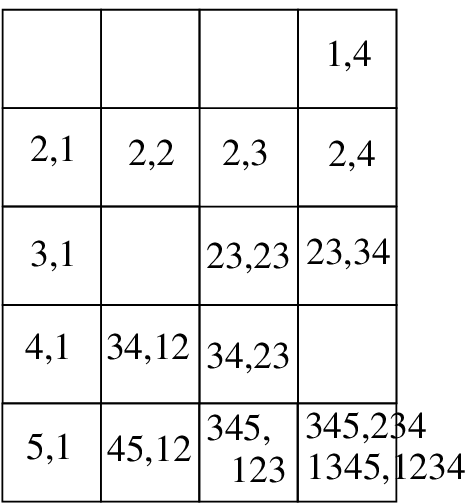
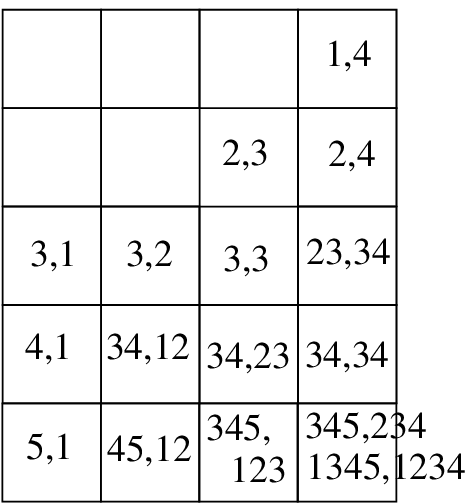
\caption{(T9),(T10),(T11) }
\label{fig:letrasf2-3}
\end{figure}

(T9) Exchange $\Delta_{12,14}$ with $\Delta_{23,34}$. That is, an NE-jump from $2$-rkline to $3$-rkline.

$$\Delta_{12,14}\Delta_{23,34} = \Delta_{123,134}\Delta_{2,4}$$

Then sweep through the $3$-rkline.

(T10) Exchange $\Delta_{134,123}$ with $\Delta_{345,234}$. That is, an SW-jump from $3$-rkline to $4$-rkline. We have obtained a $J\LE$-diagram.

$$\Delta_{134,123}\Delta_{345,234} = \Delta_{1345,1234}\Delta_{34,23}$$

(T11) Fold at $\Delta_{34,23}$, fold at $\Delta_{2,1}$, $\Delta_{2,2},\Delta_{23,12},\Delta_{23,23}$ to get a $G\LE$-diagram.

$$\Delta_{34,23}\Delta_{23,12} = \Delta_{34,12}\Delta_{23,23}$$
$$\Delta_{2,1}\Delta_{3,2} = \Delta_{2,2}\Delta_{3,1} + \Delta_{23,12}\Delta_{\phi,\phi}$$
$$\Delta_{2,2}\Delta_{3,3} = \Delta_{2,3}\Delta_{3,2} + \Delta_{23,23}\Delta_{\phi,\phi}$$
$$\Delta_{23,12}\Delta_{34,23} = \Delta_{34,12}\Delta_{23,23}$$
$$\Delta_{23,23}\Delta_{34,34} = \Delta_{23,34}\Delta_{34,23}$$

So we want a transformation process $D \rightarrow SD \rightarrow JD \rightarrow GD$ where
\begin{itemize}
\item $D$ is a $\LE$-diagram of a cell.
\item $SD$ is a $S\LE$-diagram of the cell. This is obtained by applying straightening process to $D$.
\item $JD$ is a $J\LE$-diagram of the cell. This is obtained by applying refining and jumping process to $SD$.
\item $GD$ is a $G\LE$-diagram of the cell. This will look like a $\LE$-diagram obtained by deleting a NW-corner of $1$-rkline of $D$, except for one unique-path-varible of $D$.

\end{itemize}

\section{Proof of Lemma~\ref{lem:main} II : Transforming $\LE$-diagram to a $S\LE$-diagram}

Now let's describe each step of the algorithm explicitly. Fix any weakly-connected cell inside $Gr_{rn}^{tnn}$. Denote the $\LE$-diagram $D_1$. Pick an uppermost NW-corner of the $1$-rkline and denote it $\alpha=(a,b)$. Let $m_{\alpha}$ denote the maximal rank among dots dominated by $\alpha$. 

Denote the $\LE$-diagram obtained by deleting $\alpha$ as $D_2$. Recall that we want to transform $D_1$ to a TP-diagram that looks like $D_2$ with an extra rank $m_{\alpha}$-variable. So we have to replace all variables of form $\Delta_{aX,bY}$ with $\Delta_{X,Y}$. And instead of $\Delta_{a,b}$, we should have variable of form $\Delta_{aA,bB}$ which corresponds to a unique path. In this section, we will show an algorithm to go from a $\LE$-diagram to a $S\LE$-diagram. This will be called the straightening process.

\begin{remark}
Throughout the paper if we state something like \textit{for SW(NE)-object a property holds}, it means that the property holds for SW-object and NE-object where SW stands for southwest and NE stands for northeast. We will only prove the statement for SW-object in such case, since the same argument can be used to prove for NE-objects.
\end{remark}

During the straightening process, we will only use the folding procedure. The resulting TP-diagram after our process will look like the following. Pick the NE and SW-end of the $k/\alpha$-rkline. Let them be at $(x_1,y_1),(x_2,y_2), x_1 > x_2, y_1 < y_2$. Then draw a line by going from $(x_1,y_1)$ to $(x_2,y_1)$ then to $(x_2,y_2)$. Then after the straightening process, all $k$-rank variables dominated by $\alpha$ will be on this line. That is, this line will be our $k/\alpha$-rkline after the 
straightening process.

\begin{lemma}
\label{lem:endrklineadjrk}
Pick the SW(NE)-end of the $k/\alpha$-rkline in $D_1$, name it $\beta$. If there is a adjacent dot to the left(above), it has rank $k-1$.
\end{lemma}
\begin{proof}
Call the adjacent dot to the left of $\beta$ as $\gamma$. Then $\gamma$ is either dominated by $\alpha$ or dot below $\alpha$. So we can ignore the region left of $\alpha$. Assume $\gamma$ is of rank $k$. Since $\gamma$ is not in $k/\alpha$-rkline, it is dominated by dot below $\alpha$. Denote $\gamma$ by $\Delta_{\{x_1,\cdots,x_k\},\{y_1,\cdots,y_k\}}$. Then by stability condition, there exits dots at $(x_1,y_1),\cdots,(x_k,y_k), x_1 < a$. Now denote the dot of rank $k$ dominating $\beta$ as $\beta_k$. Then by $\LE$-condition, $\beta_{k-1}$ cannot lie in region right of $(x_k,y_k)$. And it cannot lie in region left of $(x_k,y_k)$ or otherwise $\beta_1=\alpha$ will dominate $\gamma$. From the fact that $\beta_k$ is right of $(x_k,y_k)$, each $\beta_t$ has to lie above $(x_{t+1},y_{t+1})$ for $1 \leq t \leq k-1$. But this implies $\beta_1=\alpha$ is above $(x_2,y_2)$, a contradiction.
\end{proof}

\begin{lemma}
\label{lem:foldableempty}
Pick a smooth TP-diagram. Pick a SE-corner $\beta$ of $k/\alpha$-rkline. If it does not have adjacent dots to the below or right, then we can fold at $\beta$. 
\end{lemma}

\begin{proof}
Look at the Pl\"ucker mutation corresponding to a fold. Let's denote $\beta$ by $\Delta_{X,Y}$ and adjacent dots of the SE-corner as $\Delta_{X',Y}$ and $\Delta_{X,Y'}$. Due to the condition given, $\Delta_{X \cup X',Y \cup Y'}$ is zero. So we can fold at $\beta$.
\end{proof}

 Pick any SE-corner as given in the Lemma~\ref{lem:foldableempty}. The folded $k$-rkline may touch the $(k-1)$-rkline but would never overlap due to the fact that we started out from $\LE$-diagram. Now if it does touch the $(k-1)$-rkline, we fold the $(k-1)$-rkline at that point due to the following lemma.
 
\begin{lemma}
\label{lem:foldablecon}
If a folded $k$-rkline touches a $(k-1)$-rkline, it is at a SE-corner of the $(k-1)/\alpha$-rkline. And we can fold again at that SE-corner of the $(k-1)$-rkline.
\end{lemma}
\begin{proof}
Follows from the pucker mutation corresponding to the fold. The folded dot of the $k/\alpha$-rkline provides us with $\Delta_{X \cup X',Y \cup Y'}$. The dot dominating the dot of SE-corner of $(k-1)/\alpha$-rkline gives us $\Delta_{X \cap X',Y \cap Y'}$. And two adjacent dots of the SE-corner of $(k-1)/\alpha$-rkline gives us $\Delta_{X',Y}$ and $\Delta_{X,Y'}$.
\end{proof}

When we fold at dot $\beta$ and whenever the folded rkline touches another rkline, fold again as in Lemma~\ref{lem:foldablecon}, we will call this procedure the fold-chain at $\beta$.

\begin{lemma}
When we do a fold-chain on a dot that doesn't dominate any dot in the TP-diagram, we get a pre-TP-diagram that is smooth and stable. That is, we get a stable TP-diagram.
\end{lemma}
\begin{proof}
It is direct from the definitions that each fold preserves the smoothness property. Each fold may break the stability, but after we finish a fold-chain, the resulting pre-TP-diagram is stable.
\end{proof}

So we first do a fold-chain on all SE-corners of $m_{\alpha}/\alpha$-rkline, which is possible due to Lemma~\ref{lem:foldableempty}. Assume there is a $k$ such that $t/\alpha$-rkline for $k+1 \leq t \leq m_{\alpha}$ has a unique NW-corner (We will call this that $t/\alpha$-rkline is straightened) but $k/\alpha$-rkline is not. Then we can again apply fold-chains on SE-corners of $k/\alpha$-rkline by the following lemma.

\begin{lemma}
When for $k+1 \leq t \leq m_{\alpha}$, $t/\alpha$-rklines are straightened, we can start a fold chain at SE-corners of the $k/\alpha$-rkline. 
\end{lemma}
\begin{proof}
Denote position of each $\alpha_t$ as $(a_t,b_t)$. Pick a SE-corner of $k/\alpha$-rkline $\beta$. Denote it $\Delta_{Xx,Yy}$.  Either $x < a_t$ or $y < b_t$. Assume $y < b_t$. Denote the adjacent dots on $k/\alpha$-rkline as $\Delta_{X'x',Yy},\Delta_{Xx,Y'y'}$. We will prove $\Delta_{X'x'x,Y'y'y}$ is zero. 

If $\Delta_{X'x'x,Y'y'y}$ is nonzero, we have non-vertice-crossing paths from $X' \cup \{x',x\}$ to $Y' \cup \{y',y\}$. The path going from $x$ to $y$ should not cross path going from $x'$ to $y'$. This implies that in the $\LE$-diagram, we have a dot at $(z,y), z  \geq a_t$ not dominated by $\alpha$. This dot is either rank $k$ or $k+1$. But it cannot of rank $k$ since it is below $\beta$ which is a SE-corner of $k$-rkline. It cannot be of rank $k+1$ since if it is, it is adjacent to the SW-end of $(k+1)/\alpha$-rkline. But this contradicts Lemma~\ref{lem:endrklineadjrk}. So $\Delta_{X'x'x,Y'y'y}$ is zero and we can fold at $\beta$.

Stability still holds for all dots of $t$-rklines where $1 \leq t \leq k-1$, we can use Lemma~\ref{lem:foldablecon}. So we can start a fold chain at $\beta$.

\end{proof}

So we apply fold chains on SE-corner of $k/\alpha$-rkline, then repeat the above process. Then all rklines will eventually be straightened.

\begin{lemma}
\label{lem:slehcond}
Applying numerous fold-chains, we get a TP-diagram such that each $k/\alpha$-rkline has no SE-corner. Denote such diagram obtained from $D_1$ as $SD_1$. Then dot $(x,y)$ is in $k/\alpha$-rkline of $SD_1$ if and only if there were $k$-rank dots $(x',y),(x,y')$ in $k/\alpha$-rkline of $SD_1$ with $x' \geq x, y' \geq y$.
\end{lemma}
\begin{proof}
Follows directly from process above.
\end{proof}

Let's call the TP-diagram obtained from $\LE$-diagrams after the straightening process as $S\LE$-diagrams. S stands for straightened.

\begin{remark}
\label{rem:nochange}
The dots not dominated by $\alpha$ and SW(NE)-ends of $k/\alpha$-rklines have not moved in the S$\LE$-diagram, compared to the $\LE$-diagram.
\end{remark}

Let's introduce some additional notation we will be using often. For TP-diagrams that each $k/\alpha$-rkline has unique NW-corners for all rank $k$, we will be using the following notation.

\begin{figure}[ht]
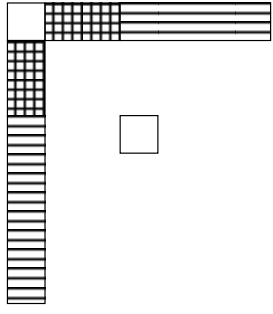
\caption{$\alpha_k,\alpha_{k+1}$ are represented by white boxes. $k$-SW(NE)-lblock is represented by box filled with horizontal lines. $k$-SW(NE)-sblock is represented by box filled with squares. }
\label{fig:block}
\end{figure}

\begin{itemize}
\item $\alpha_k$ stands for unique NW-corners of $k/\alpha$-rklines.
\item $(a_k,b_k)$ denotes position of $\alpha_k$.
\item $k$-block stands for $k/\alpha$-rkline.
\item $k$-SW(NE)-block stands for dots in $k/\alpha$-rkline below(right) of $\alpha_k$.
\item $k$-SW(NE)-lblock stands for dots in $k$-SW(NE)-block that are in rows(columns) weakly below(right) of row(column) $\alpha_{k+1}$ is in.
\item $k$-SW(NE)-sblock stands for complement of $k$-SW(NE)-lblock inside $k$-SW(NE)-block.

\end{itemize}

\begin{lemma}
\label{lem:semistability}
If SW(NE)-ends of $k$-block and $(k+1)$-block are connected, then $k$-NE(SW)-sblock is empty. And if $k$-NE(SW)-sblock is empty, the following holds. Let $\beta$ denote either $\alpha_{k+1}$ or a dot of $(k+1)$-SW(NE)-block. Pick a lowermost(rightmost) dot $\gamma$ among dots of $k$-SW(NE)-block that lie in region above(left) of $\beta$. Then $\gamma$ dominates $\beta$.
\end{lemma}

\begin{proof}
First statement is trivial. For the second statement, since $(k+1)$-NE(SW)-sblock is empty, $\gamma$ was folded only when $\beta$ was folded. So it means $\gamma$ dominated $\beta$ in the $\LE$-diagram.
\end{proof}

\begin{lemma}
For each dot on $k$-SW(NE)-block, draw a line to the right(below). Then there are dots on the $k'$-blocks, where $k'>k$, on the intersection with that horizontal(vertical) line. \end{lemma}

\begin{proof}
By Lemma~\ref{lem:slehcond} and the folding process.
\end{proof}

\section{Proof of Lemma~\ref{lem:main} III : Transforming $S\LE$-diagram to a $J\LE$-diagram}

We follow the notations used in the previous section. In the previous section, we have obtained an algorithm to transform a $\LE$-diagram to a $S\LE$-diagram. In this section, we will obtain an algorithm to transform it to a $J\LE$-diagram. The process involved will be called the refining process and jumping/sweeping process.

Let's look at the refining process. Recall that we have fixed $D_1$ and obtained $SD_1$ from it. Let's look at some properties of $SD_1$.

\begin{definition}

Pick any two adjacent points of $k$-SW-block. Let them be at position $(x,y),(x',y)$, $x < x'$ and variables be $\Delta_{aXx,bYy},\Delta_{aX'x',bYy}$. 
We will say $\Delta_{aX'x',bYy}$ is unswappable if $\Delta_{Xxx',bYy}  >0$ and swappable if it is zero. We will say $\Delta_{aX'x',bYy}$ is uploadable if $\Delta_{Xx,bY}> 0$. So uploadable means the following. If one of the following conditions are satisfied,
\begin{itemize}
\item $\Delta_{X'x,bY}$ or $\Delta_{Xxx',bYy}$ is zero in the cell.
\item $\Delta_{X'x,bY}$ and $\Delta_{Xxx',bYy}$ is inside the diagram.
\end{itemize}
we can push at $\Delta_{aX'x',bYy}$. Swappability and uploadability is defined similarly for two adjacent vertices of $k$-NE-block.

\end{definition}

\begin{lemma}
If a dot is unswappable, it is uploadable.
\end{lemma}
\begin{proof}
Follows directly from the definitions.
\end{proof}

\begin{lemma}
\label{lem:uploadability}
A point in a $k$-SW(NE)-block is uploadable if and only if SW(NE)-end of $k$-SW(NE)-block is adjacent to a $k-1$ rank dot not dominated by $\alpha$.
\end{lemma}

\begin{proof}
Denote $\beta$ by $\Delta_{aX'x',bYy}$. Let its uploadability condition be $\Delta_{Xx,bY} >0$. Let $Xx = \{x_1,\cdots,x_{k-1}\}, bY = \{y_1,\cdots,y_{k-1}\}$ such that $x_1 < \cdots < x_{k-1}, y_1< \cdots < y_{k-1}$. From $\Delta_{Xx,bY}>0$, we have dots $(w_1,y_1),\cdots,(w_{k-1},y_{k-1})$ such that $w_i \geq x_i$ for $1 \leq i \leq k-1$ in $D_1$. Denote dot at $(w_{k-1},y_{k-1})$ by $\gamma$. Then $\gamma$ is not dominated by $\alpha$. So $\gamma$ has to lie on row weakly below row of SW-end of $k$-block. But if it is strictly below the row of SW-end of $k$-block, we don't have a path going from $x_{k-1}$ to $y_{k-1}$. So $\gamma$ is adjacent to SW-end of $k$-block. By Lemma~\ref{lem:endrklineadjrk}, $\gamma$ is of rank $k-1$. This proves the only if direction. By Remark~\ref{rem:nochange} and the supersmoothness of $\LE$-diagram, $\gamma$ is $\Delta_{H \setminus \{a\},J \setminus max(J)}$ where SW-end of SW-block is denoted as $\Delta_{H,J}$. 

Now denote the SW-end of $k$-SW-block as $\Delta_{aX_1,Yy}$. Let the dots above that on the $k$-SW-block be $\Delta_{aX_2,Yy},\Delta_{aX_3,Yy},\cdots$ from bottom to top. The $k-1$-rank dot to the left of $\Delta_{aX_1,Yy}$ is $\Delta_{X_1,Y_1}$ due to the above argument.
$$\Delta_{aX_1,Yy} \Delta_{X_2,Y} =  \Delta_{X_1,Y} \Delta_{aX_2,Y,y} + \Delta_{X_1 \cup X_2,Yy} \Delta_{\{a, X_1 \cap X_2\},Y}$$
gives us $\Delta_{X_2,Y} >0$. So $\Delta_{aX_1,Yy}$ is uploadable. Now if we look at the following relation for each $k$,
$$\Delta_{aX_k,Yy} \Delta_{X_{k+1},Y} =  \Delta_{X_k,Y} \Delta_{aX_{k+1},Y,y} + \Delta_{X_k \cup X_{k+1},Yy} \Delta_{\{a, X_k \cap X_2\},Y}$$
tells us each $\Delta_{aX_k,Yy}$ is uploadable. So all dots of $k$-SW-block are uploadable. 
\end{proof}

\begin{corollary}
Pick any point on $k$-SW(NE)-block. If that point is uploadable, then all points in that block is uploadable.
\end{corollary}

By above corollary, we can say a $k$-SW(NE)-block is uploadable or non-uploadable, since all dots in the $k$-SW(NE)-block are uploadable or non-uploadable. 

\begin{corollary}
\label{cor:sblockempty}
If $(k+1)$-SW(NE)-block is uploadable, then $k$-NE(SW)-lblock is empty.
\end{corollary}
\begin{proof}
Follows from Lemma~\ref{lem:uploadability}, $\LE$-condition and Remark~\ref{rem:nochange}. \end{proof}

So by Lemma~\ref{lem:semistability}, if $(k+1)$-SW(NE)-block is uploadable, each dot in $(k+1)$-SW(NE)-block is dominated by a dot in $k$-SW(NE)-block.

\begin{corollary}
If there is an unswappable dot in $k$-SW(NE)-block, denote it $\beta$. Let its unswappable condition be $\Delta_{H,J} >0$. This happens if and only if there is a dot in $(k+1)$-block, having its uploadability condition as $\Delta_{H,J} >0$. Call this dot helper of $\beta$. Then $\beta$ is dominates the helper of $\beta$.
\end{corollary}


So we can say a lblock is unswappable or swappable, since either all dots in the lblock are swappable or unswappable. If it is empty, we will also say that it is swappable. All dots in $k$-SW(NE)-sblocks are swappable by the above corollary. Now we want to know when an uploadable dot can be pushed. For convenience, we will say that a $k$-SW(NE)-block is unswappable if $k$-SW(NE)-lblock is unswappable.

\begin{proposition}
\label{prop:helpers}
Denote an uploadable dot in $k$-SW(NE)-block as $\beta$. If there is an adjacent dot to the below(right) on the block, call it $\gamma$. If $\beta$ is unswappable, call the helper of $\beta$ as $\omega$. Then $\beta$ can be pushed if one of the following is satisfied.

\begin{enumerate}
\item $\gamma$ and $\omega$ does not exist.
\item $\gamma$ does not exists, $\omega$ exists and is pushed.
\item $\gamma$ exists and is pushed, $\omega$ does not exists.
\item $\gamma$ and $\omega$ exist and both of them are pushed.
\end{enumerate}

\end{proposition}
\begin{proof}
Follows from the pucker mutation corresponding to the push. Denote $\beta$ by $\Delta_{aX,Yb}$ and adjacent dot to the above by $\Delta_{aX',Yb}$. Non-existence of $\gamma$ means $\beta$ is SW-end of the block. By Lemma~\ref{lem:uploadability}, the adjacent dot of $\gamma$ with rank $k-1$ gives us $\Delta_{X,Y}$. If $\gamma$ exists, $\gamma$ being pushed gives us $\Delta_{X,Y}$. Non-existence of $\omega$ gives us $\Delta_{X \cup X',Yb}=0$ in the cell. If $\omega$ exists, $\omega$ being pushed provides us with $\Delta_{X \cup X', Yb}$. $\Delta_{\{a\} \cup (X \cap X'), Y}$ is given by dot dominating $\beta$.
\end{proof}

So we can push all dots of a $k$SW(NE)-block if it is uploadable. We do this by pushing the $k$-SW(NE)-end first then moving along the block as we push the adjacent dots. It is the same for $k$-SW(NE)-lblock. We can push the $k$-SW(NE)-sblock after we have pushed the $k$-SW(NE)-lblock. Pulling is done by reversing the process of pushing.

\begin{lemma}
\label{lem:pullanytime}
It is always possible to pull back blocks that were swappable but uploaded.
\end{lemma}
\begin{proof}
Look at the Pl\"ucker mutation of the pull. Since in this case, $\Delta_{X \cup X', Yb}$ is zero due to the fact that we are looking at a block that is swappable. So we can always push or pull the block.
\end{proof}

Our refining process will be pushing necessary blocks to make the jumping from $\alpha$ to all other dots dominated by $\alpha$ possible. We push all blocks that are unswappable. There are blocks that are swappable but uploadable. We want to decide whether or not to push them.

\begin{definition}
Let us have $\Delta_{aX,bY},\Delta_{X,Y},\Delta_{aX',bY'}$ in the pre-TP-diagram but not $\Delta_{X',Y'}$. We will say that we jump from dot $\Delta_{aX,bY}$ to $\Delta_{aX',bY}$ by replacing $\Delta_{aX,bY}$ with $\Delta_{X',Y'}$.
\end{definition}

So what we want to do is the following. Starting from $\alpha$, we want a sequence of jumps that goes through all dots dominated by $\alpha$ of the pre-TP-diagram, not landing on the same dot twice. Sweeping is jumping through adjacent dots of a block. When all dots dominated by $\alpha$ except $\alpha_k$'s are swappable, we can sweep through the $k$-blocks in any direction(From SW-end to NE-end or from NE-end to SW-end). We are going to prove later that when SW(NE)-ends of $k$-SW(NE)-block and $(k+1)$-SW(NE)-block of swappable blocks are connected, we can jump from SW(NE)-end of $k$-SW(NE)-block to $(k+1)$-SW(NE)-block. We are looking at weakly-connected cells. Assume the following conditions hold.

\begin{itemize}
\item All blocks are swappable.
\item For $k$ odd, SW-ends of $k$-block and $(k+1)$-block are connected.
\item For $k$ even, NE-ends of $k$-block and $(k+1)$-block are connected.
\end{itemize}

Then we could process as the following.

\begin{enumerate}
\item Jump from $\alpha$ to SW-end of $2$-block.
\item Sweep from SW-end to NE-end of $2$-block.
\item Jump from NE-end of $2$-block to NE-end of $3$-block.
\item Sweep from NE-end to SW-end of $3$-block.
\item Jump from SW-end of $3$-block to SW-end of $4$-block.
\item ... So on.
\end{enumerate}

But the first condition we assumed normally does not hold. We can't sweep in any direction of $k$-block if a $k$-SW(NE)-lblock is unswappable. This is why we push unswappable blocks. 

\medskip

Recall that $m_{\alpha}$ denotes the maximal rank among all dots dominated by $\alpha$ in $\LE$-diagram of the cell. Since we have fixed $\alpha$, denote $m_{\alpha}$ by $m$. If $m=1$, there is nothing to do, so assume $m \geq 2$.

Now let's start the refining process. Denote $k_1,k_2$ to be the minimal rank of SW,NE blocks which are swappable. Then SW(NE)-$k_1(k_2)$-block is swappable and uploadable. 
Assume $k_1 \leq k_2$. So SW-blocks of rank $2,\cdots,k_1$ are uploadable and NE-blocks of rank $2,\cdots,k_2$ are uploadable. Push $(k_1)$-SW-block, then $(k_1-1)$-SW-block, and repeat till we push $2$-SW-block. This is possible due to Proposition~\ref{prop:helpers}. Then pull back $k_1$-SW-block which is possible due to Lemma~\ref{lem:pullanytime}. Do the same for NE-blocks.

The situation can more simply be described as a graph on a grid. To avoid confusion, we denote positions of points in this grid as $<x,y>$. When we delete a dot, we also delete the adjacent edges. 

\begin{figure}[ht]
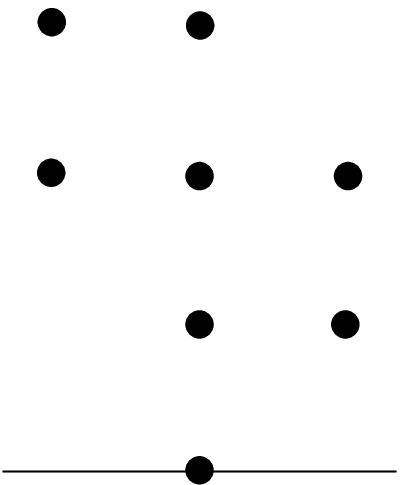
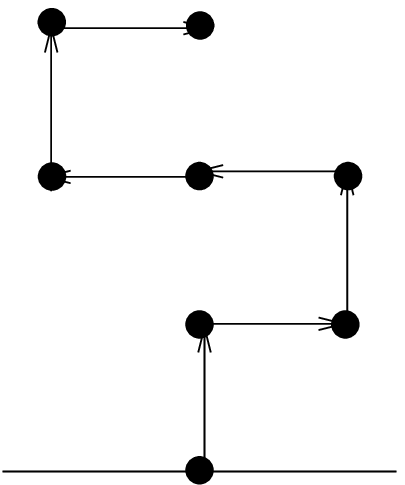
\caption{Second example in Section 7 and its Hamiltonian path.}
\label{fig:3dotex}
\end{figure}

\begin{enumerate}
\item On a grid, draw dots at $<0,1>,<0,2>, \cdots, <0,m>$. They stand for $\alpha,\cdots,\alpha_m$.
\item Now draw dots $<-1,k_1>,\cdots,<-1,m>,<1,k_2>,\cdots,<1,m>$. They stand for SW,NE-ends of each block. Note that when $k$-SW(NE)-block was empty, then SW-end of $k$-block is also $\alpha_k$. So $<-1,k>$ and $<0,k>$ both denotes $\alpha_k$. $\alpha_4$ of second example in Section 7 is such an example. See Figure~\ref{fig:3dotex}.

\item If $k_1=2$, draw edge $<0,1> \rightarrow <-1,2>$. If $k_2=2$, draw edge $<0,1> \rightarrow <1,2>$. This corresponds to the fact $\alpha_1$ is SW-end and NE-end of $1$-block.

\item For $k \geq k_1$, draw $<-1,k> \rightarrow <-1,k+1>$ if SW-ends of $k,(k+1)$-block are connected. The edges drawn correspond to the fact that we can jump from SW-end of $k$-block to SW-end of $(k+1)$-block for $k \geq k_1$, which will be proven in Lemma~\ref{lem:endjump}. 
\item Similarly, for $k \geq k_2$, draw $<1,k> \rightarrow <1,k+1>$ if NE-ends of $k,(k+1)$-block are connected. 
\item Now draw edges $<0,k> \rightarrow <0,k+1>$ for $1 \leq k < k_1-1$. This corresponds to fact that when SW and NE blocks of $(k+1)$-block are pushed, it is jumpable from $\alpha_k$ to $\alpha_{k+1}$. This will be proven in Corollary~\ref{cor:midjump}. 
\item Weakly-connectedness tells us we whenver we have vertices $< (-1)^{k+t} , k>, < (-1)^{k+t},k+1>$ we have edge  $< (-1)^{k+t} , k> \rightarrow < (-1)^{k+t},k+1>$ where $t$ is either $0$ or $1$.
\item Orientation of the horizontal edges will correspond to which direction we do the sweeping of each $k$-blocks. We have to decide them
\end{enumerate}

 Define $par(x)$ as the parity of integer $x$. For $k > k_2$ we can either draw $<-1,k>\rightarrow <0,k>\rightarrow  <1,k>$ or $<-1,k> \leftarrow <0,k> \leftarrow <1,k>$ due to the fact that $k$-block is swappable.
\begin{itemize}
\item If $par(k) = par(t)$, draw edges $<-1,k> \rightarrow <0,k>, <0,k> \rightarrow <1,k>$. 
\item If $par(k) \not = par(t)$, draw edges $<1,k> \rightarrow <0,k>, <0,k> \rightarrow <-1,k>$.

\end{itemize}

Now we have to decide whether to push SW(NE)-$k_1$($k_2$)-block or not. We express them by $<-1,k_1>$ and $<1,k_2>$. $<-1,k_1>$ corresponds to $k_1$-SW-block. If we push the $k_1$-SW-block, we delete the dot. $<1,k_2>$ corresponds to $k_2$-NE-block. If we push the $k_2$-NE-block, we delete the dot.

\begin{itemize}

\item If $k_1 = k_2$.
\begin{itemize}
\item If $par(k_1) = par(t)$. We have $<1,k_1> \rightarrow <1,k_1+1>$.
Push $k_1$-SW-lblock and delete $<-1,k_1>$. Draw edge $<0,k_1-1> \rightarrow <0,k_1>$ which corresponds to the fact that we can do a sequence of jumps from $\alpha_{k_1}$ to SW-end of $k_1$-SW-sblock, then to $\alpha_{k+1}$ due to Lemma~\ref{lem:midjump}.

\item If $par(k_1) \neq par(t)$. We have $<-1,k_1> \rightarrow <-1,k_1+1>$. Push $k_1$-NE-lblock and delete $<1,k_1>$. Draw edge $<0,k_1-1> \rightarrow <0,k_1>$ which corresponds to the fact that we can do a sequence of jumps from $\alpha_{k_1}$ to NE-end of $k_1$-NE-sblock, then to $\alpha_{k+1}$ due to Lemma~\ref{lem:midjump}.
\end{itemize}

\item If $k_1 < k_2$.

For $<-1,k_1>$,

\begin{itemize}
\item If $par(k_1) \neq par(t)$. We have $<-1,k_1> \rightarrow <-1,k_1+1>$. Do not push the $k_1$-SW-block. Draw edge $<-1,k_1> \leftarrow <0,k_1>$.
\item If $par(k_1) = par(t)$. Push the $k_1$-SW-block and delete $<-1,k_1>$.
 
\end{itemize}

For $k_1 < k < k_2$, draw edge 

\begin{itemize}
\item $<-1,k> \rightarrow <0,k>$ if $par(k) = par(t)$
\item $<-1,k> \leftarrow <0,k>$ if $par(k) \neq par(t)$.
\end{itemize}

For $k_1-1 \leq k \leq k_2-1$, draw edge $<0,k> \rightarrow <0,k+1>$ if $par(k)=par(t)$. This corresponds to the following. Pull the $k$-NE-sblock, which is possible due to Lemma~\ref{lem:pullanytime}. Do a sequence of jumps from $\alpha_k$ to NE-end of $k$-NE-sblock. Then we can jump from that dot to $\alpha_{k+1}$ by Lemma~\ref{lem:midjump}. 
 
For $<1,k_2>$,
\begin{itemize}
\item If $par(k_2) = par(t)$. We have $<1,k_2> \rightarrow <1,k_2+1>$. Do not push the $k_2$-NE-block. Draw edges $<-1,k_2> \rightarrow <0,k_2> \rightarrow <1,k_2>$.
\item If $par(k_2) \neq par(t)$. Push the $k_2$-NE-block and delete $<1,k_2>$.
\end{itemize}

\end{itemize}

This finishes the refining process, and we get the transformed diagram and graph $G'$ on a grid. This $G'$ will act as a map for our jumping process. From above construction, $G'$ has a Hamiltonian path $P'$ starting from $<0,1>$ and ending at last row, going through all vertices. The jumping goes as the follows. Start from $\alpha$, set $i=2$, start the following procedure.

If there is a non-horizontal edge in $P'$ heading to
\begin{itemize}
\item $<0,i>$, jump to $\alpha_i$.
\item $<-1,i>$, jump to SW-end of $i$-block.
\item $<1,i>$, jump to NE-end of $i$-block.
\end{itemize}

If $i$-th row of $G'$ is oriented to the right(left), sweep in NE(SW)-direction of the current block. Increase $i$ by $1$, repeat the above procedure until we have reached the last dot dominated by $\alpha$. So using the $P'$ as a map, we have finished the jumping process.

 Let's see an example of the graph constructed above. Look at the cell selected at second example of Section 7. The $2$-SW-block was unswappable, so we pushed it. We did not push $3$-SW-block and $2$-NE-block, so we didn't delete $<-1,3>$ and $<1,2>$.

 So we have obtained a TP-diagram such that only one dot is of form $\Delta_{aX,bY}$. Denote it $\beta$. $\beta$ is a unique-path-variable because it is of rank $m$. Call this the $J\LE$-diagram of the cell. Denote this $J\LE$-diagram obtained from $SD_1$ as $JD_1$. We want to transform this $J\LE$-diagram to look like $D_2$, except for $\beta$. We will do this in the next section.
 
 But before we go on, we have to prove Lemma~\ref{lem:midjump},Lemma~\ref{lem:endjump} and Corollary~\ref{cor:midjump}.

\begin{lemma}
If both(SW and NE) sblocks of $k$-block are empty, then it is jumpable from $\alpha_k$ to $\alpha_{k+1}$.
\end{lemma}

\begin{proof}
If both sblocks of $k$-block are empty, this implies $\alpha_{k+1}$ is dominated by $\alpha_k$. 
In notation of the Pl\"ucker mutations of the fold, denote $\alpha_k = \Delta_{X,Y}, \alpha_{k+1} = \Delta_{aX',bY'}$. To jump from $\alpha_k$ to $\alpha_{k+1}$, we need to exchange $\Delta_{X,Y}$ with $\Delta_{X',Y'}$.
$\alpha_{k+1}=\Delta_{X \cup X',Y \cup Y'}$. And $\Delta_{X \cap X',Y \cap Y'}$ is given since we are assuming we are jumping from $\alpha_{k}$. So to exchange $\Delta_{X,Y}$ with $\Delta_{X',Y'}$, we need $\Delta_{X',Y}$ and $\Delta_{X,Y'}$ if both of them are nonzero. 

Look at $k$-SW-block and $k$-NE-block.

\begin{itemize}
\item If both of them are pushed, the dots that were adjacent to $\alpha_{k+1}$ before the pushing provides $\Delta_{X',Y}$ and $\Delta_{X,Y'}$.

\item If one of them is not uploadable, then $\Delta_{X',Y}$ or $\Delta_{X,Y'}$ is zero in the cell.

\item If $k$-SW-block is uploadable but not pushed, this means this block is swappable. So we can temporarily upload it to get $\Delta_{X,Y'}$. After we have jumped, pull the block back again, which is possible due to Lemma~\ref{lem:pullanytime}. This is exactly what we did in (T3),(T5) of both examples in Section 7. Similar for $k$-NE-block.
  
\end{itemize}

So we can do the fold to jump from $\alpha_k$ to $\alpha_{k+1}$.

\end{proof}

From above lemma, we get the following corollary.

\begin{corollary}
\label{cor:midjump}
If both $(k+1)$-SW-block and $(k+1)$-NE-block are uploadable, it is jumpable from $\alpha_k$ to $\alpha_{k+1}$
\end{corollary}

\begin{lemma}
\label{lem:midjump}
If $k$-NE(SW)-sblock is empty but $k$-SW(NE)-sblock is not, it is jumpable from SW(NE)-end of $k$-SW(NE)-sblock to $\alpha_{k+1}$. 
\end{lemma}
\begin{proof}
If NE-sblock of $k$-block is empty but SW-sblock is not, this implies SW-end of $k$-SW-sblock dominates $\alpha_{k+1}$. Denote SW-end of $k$-SW-sblock as $\Delta_{aX,bY}$ and $\alpha_{k+1}$ as $\Delta_{aXa_k,bYb_k}$. We want to prove $\Delta_{aX,Yb_k}=0$. 

Assume $\Delta_{aX,Yb_k} > 0$. 

Denote $aX = \{w_1,\cdots,w_k\},Yb_k = \{z_1,\cdots,z_k\}$ such that $w_1 < \cdots < w_k, z_1 < \cdots < z_k$. Then we have dots in $D_1$ at $(v_1,z_1),\cdots,(v_k,z_k)$ where $v_i \leq w_i$. By $\LE$-condition we have a dots at $(v_{t-1},z_t)$, $2 \leq t \leq k$. Since $\alpha_k$ is at column $z_{t-1}$, by smoothness implies $(v_{t-1},z_t)$ has rank $k+1$. So pick the uppermost dot of rank $k+1$ on column $z_k$. This dot is dominated by a dot of $k$-block.

So $\Delta_{aX,Yb_k}>0$ implies there is a rank $k+1$ dot dominated by $\alpha$ on column $max(Y)$. This contradicts the fact that $\alpha_{k+1}$ is unique NW-corner of $(k+1)$-block. Hence $\Delta_{aX,Yb_k}=0$. Now looking at the Pl\"ucker relation corresponding to a fold, we have the following relation.
$$\Delta_{Xa_k,Yb_k}\Delta_{aX,bY} = \Delta_{aXa_k,bYb_k}\Delta_{X,Y}$$

Since we have $\alpha_{k+1}=\Delta_{aXa_k,bYb_k}, \Delta_{X,Y}$ in the current TP-diagram, we can exchange $\Delta_{aX,bY}$ with $\Delta_{Xa_k,Yb_k}$. 
\end{proof}

\begin{lemma}
\label{lem:endjump}
If SW(NE)-ends of $k$-block and $(k+1)$-block are connected, and $k$-SW(NE)-block is swappable, we can jump from SW-end of $k$-block to SW(NE)-end of $(k+1)$-block.
\end{lemma}
\begin{proof}
Recall that we called SW-ends of $k,(k+1)$-block connected if one of the following holds.
\begin{enumerate}
\item SW-end of $k$-block dominates SW-end of $(k+1)$-block
\item SW-end of $k$-block is adjacent to SW-end of $(k+1)$-block
\end{enumerate}

$k$-SW-end being swappable implies $(k+1)$-SW-end unuploadable. 

For case (1), denote $k$-SW-end as $\Delta_{aX,bY}$ and $(k+1)$-SW-end as $\Delta_{aXx,bYy}$. $(k+1)$-SW-end unuploadable implies $\Delta_{Xx,bY}=0$. So looking at the Pl\"ucker relation corresponding to a fold, we have the following relation.
$$\Delta_{Xx,Yy}\Delta_{aX,bY} = \Delta_{aXx,bYy}\Delta_{X,Y}$$
Since we have $\Delta_{aXx,bYy}, \Delta_{X,Y}$ in the current TP-diagram, we can exchange $\Delta_{aX,bY}$ with $\Delta_{Xx,Yy}$. 

For case (2), the proof is similar.
\end{proof}

\section{Proof of Lemma~\ref{lem:main} IV : Transforming a $J\LE$-diagram to a $G\LE$-diagram }

In this section, we will finish proving Lemma~\ref{lem:main}. Recall that in the previous sections, we have fixed a cell, its corresponding $\LE$-diagram $D_1$, its NE-corner of $1$-rkline $\alpha=\Delta_{\{a\},\{b\}}$, and a $\LE$-diagram $D_2$ obtained from $D_1$ by deleting dots corresponding to $\alpha$. 

We first will show an algorithm to transform $D_2$, so that all variables are same except for one compared to $JD_1$ of the previous section. We will imitate the straightening, refining process we have done for $D_1$. 

\begin{definition}
For pre-TP-diagram $T$ containing $\alpha$, we say that pre-TP-diagram $T'$ is obtained by trimming $\alpha$ from $T$ if we get $T'$ from $T$ by deleting $\alpha$ and changing $\Delta_{aX,bY}$ to $\Delta_{X,Y}$. Since we fixed $\alpha$ in this section, we will say that $T'$ can be obtained from $T$ by trimming.
\end{definition}

Our basic strategy is as the following. Whenever we do a fold at $\Delta_{aX,bY}$ of $D_1$, we do a fold at $\Delta_{X,Y}$ of $D_2$. Whenever we do a push at $\Delta_{aX,bY}$ of $D_1$, we do a fold at $\Delta_{X,Y}$ of $D_2$. For each step, we will denote this as imitating fold(push) of $\Delta_{aX,bY}$. Let's first prove that the entire straightening process can be imitated.

\begin{lemma}
The straightening process for $D_1$ can be imitated at $D_2$. That is, if $T_1=D_1 \rightarrow T_2 \rightarrow T_3 \cdots \rightarrow T_s=SD_1$ is the sequence of pre-TP-diagrams obtained in the straightening process, we can imitate each step to get $T_1'=D_2 \rightarrow T_2' \rightarrow T_3' \cdots \rightarrow T_s'=SD_1'$. And $T_s'$ is obtained from $T_s$ by trimming.
\end{lemma}

\begin{proof}

By definition, we know $T_1'$ is obtained from $T_1$ by trimming. Assume we already have constructed $T_1' \rightarrow T_2' \cdots \rightarrow T_i'$ by imitating $T_1 \rightarrow T_2 \cdots \rightarrow T_i$ and $T_i'$ is obtained from $T_i$ by trimming. Let's prove we can imitate $T_i \rightarrow T_{i+1}$ to get $T_i' \rightarrow T_{i+1}'$ and $T_{i+1}'$ can be obtained from $T_{i+1}$ by trimming. 

The Pl\"ucker mutation of the fold $T_i \rightarrow T_{i+1}$ can be expressed as the following relation.
$$\Delta_{aX,bY}\Delta_{aX',bY'} = \Delta_{aX',bY}\Delta_{aX,bY'} + \Delta_{aX \cap aX',bY \cap bY'}\Delta_{aX \cup aX',bY \cup bY'}$$

Now the relation we will use for $T_{i}' \rightarrow T_{i+1}'$ can be expressed as the following relation.
$$\Delta_{X,Y}\Delta_{X',Y'} = \Delta_{X',Y}\Delta_{X,Y'} + \Delta_{X \cap X',Y \cap Y'}\Delta_{X \cup X',Y \cup Y'}$$
Since $T_i'$ is obtained from $T_i$ by trimming due to induction hypothesis, we have all the necessary variables to do a Pl\"ucker mutation corresponding to this relation. And it is direct from the relation that $T_{i+1}'$ can be obtained from $T_{i+1}$ by trimming. So by induction, we have proven the lemma.

\end{proof}

\begin{lemma}
Refining process for $SD_1$ can be imitated in $SD_2$. That is, if $Q_1=SD_1 \rightarrow Q_2 \rightarrow Q_3 \cdots \rightarrow Q_f=JD_1$ is the sequence of pre-TP-diagrams obtained in the refining process, we can imitate each step to get $Q_1'=SD_2 \rightarrow Q_2' \rightarrow Q_3' \cdots \rightarrow Q_f'={JD_1}'$ where each $Q_i'$ is obtained from $Q_i$ by trimming.
\end{lemma}

\begin{proof}

By definition, we know $Q_1'$ is obtained from $Q_1$ by trimming. Assume we already have constructed $Q_1' \rightarrow Q_2' \cdots \rightarrow Q_i'$ by imitating $Q_1 \rightarrow Q_2 \cdots \rightarrow Q_i$ and $Q_i'$ is obtained from $Q_i$ by trimming. Let's prove we can imitate $Q_i \rightarrow Q_{i+1}$ to get $Q_i' \rightarrow Q_{i+1}'$ and $Q_{i+1}$ can be obtained from $Q_{i+1}$ by trimming. 

The Pl\"ucker mutation of the push/pull $Q_i \rightarrow Q_{i+1}$ can be expressed as the following relation. (Assume it is a SW-push. Proof is similar for other cases.)
$$\Delta_{aXx,bYy}\Delta_{X',bY} = \Delta_{X,bY}\Delta_{aX',bYy} + \Delta_{X \cup X', bYy}\Delta_{\{a\} \cup (X \cap X'), bY}$$

Now the relation we will use for $Q_{i}' \rightarrow Q_{i+1}'$ can be expressed as the following relation.
$$\Delta_{Xx,Yy}\Delta_{X',bY} = \Delta_{X,bY}\Delta_{X',Yy} + \Delta_{X \cup X', bYy}\Delta_{X \cap X', Y}$$
Since $Q_i'$ is obtained from $Q_i$ by trimming due to induction hypothesis, we have all the necessary variables to do a Pl\"ucker mutation corresponding to this relation. And it is direct from the relation that $Q_{i+1}'$ can be obtained from $Q_{i+1}$ by trimming. So by induction, we have proven the lemma.

\end{proof}

So we have a sequence of mutations $D_2$ to $JD_1'$. And by reversing this process, we get a sequence of mutations $JD_1'$ to $D_2$. All variables in $JD_1'$ are contained in $JD_1$. 

\begin{lemma}
\label{lem:nonzeroiff}
Let $D_2$ be a $\LE$-diagram obtained from $D_1$ by deleting $\alpha=(a,b)$, the NW-corner of $1$-rkline. Let $|H|=|J|$ such that 
\begin{enumerate}
\item $min(H) \geq a, min(J) \geq b$
\item $a \not \in H$ or $b \not \in J$.
\end{enumerate}
Then $\Delta_{H,J}$ is nonzero in $S_{D_1}$ if and only if it is nonzero in $S_{D_2}$.

\end{lemma}
\begin{proof}

Look at the paths corresponding to $\Delta_{H,J}$ inside the $\LE$-graph of $D_1,D_2$. In both graphs, the path is bounded inside the region $\{(x,y)|x \geq a, y \geq b\}$. And since $a \not \in H$ or $b \not \in J$, the path does not cross $(a,b)$. So the paths are exactly same in $\LE$-graph of $D_1$ and $D_2$.
\end{proof}

During the transformation process of the previous section, all Pl\"ucker variables we have used are dominated by $\alpha$ or rank $1$ dots adjacent to $\alpha$. So we can apply the transformation process $JD_1' \rightarrow D_2$ to $JD_1$ by the above lemma. So we can transform $JD_1$ to a TP-diagram that consists all variables of $D_2$, and a variable corresponding to a unique path in $S_{D_1}$. Denote this as $GD_1$. This proves Lemma~\ref{lem:main}.

\begin{remark}
During the transformation process $D_1 \rightarrow SD_1 \rightarrow JD_1 \rightarrow GD_1$, for all the Pl\"ucker mutations we used, only Pl\"ucker variables that are involved in the transformation process are dots dominated by $\alpha$ or rank $1$ dots adjacent to $\alpha$.
\end{remark}

\section{Proof of Theorem~\ref{thm:main} and how to use it}

Now let's prove Theorem~\ref{thm:main}. We follow the notation of the previous sections. Start from $D_1$, the $\LE$-diagram of a weakly-connected positroid cell. Obtain $D_{i+1}$ by deleting the rightmost NW-corner of $1$-rkline of $D_i$. Then if we denote the number of dots in $D_1$ as $f$, we get a sequence $D_1, D_2, \cdots, D_f$ where $D_f$ is a $\LE$-diagram with no dots inside. 

\begin{lemma}
\label{lem:multiiff}
Denote $\alpha=(a,b)$ as the rightmost NW-corner of $1$-rkline of $D_i$.  Let $|H|=|J|$ such that 
\begin{enumerate}
\item $min(H) \geq a, min(J) \geq b$
\item $a \not \in H$ or $b \not \in J$.
\end{enumerate}
Then $\Delta_{H,J}$ is nonzero in $S_{D_1}$ if and only if it is nonzero in $S_{D_i}$.
\end{lemma}

\begin{proof}
By Lemma~\ref{lem:nonzeroiff}.
\end{proof}

By Lemma~\ref{lem:main} we have transformation process $t_i : D_i \rightarrow GD_i$ where $GD_i$ is the $G\LE$-diagram of $D_i$. By Lemma~\ref{lem:multiiff} we can do the following. Apply $t_1$ to $D_1$, then $t_2$, then $t_3,\cdots,t_{f-1}$. We obtain a TP-diagram such that all of its variables are unique paths in $D_1$. So we have proven Theorem~\ref{thm:main}.

Now let's use what we have obtained to show an example of how to go from canonical TP-bases $\S$ to TP-bases $\S'$ such that all edge weights of $\LE$-network can be expressed as Laurent monomial in variables of $\S'$.

Starting from $D$ the $\LE$-diagram of the cell and $\S' = \{\}$, proceed with the following algorithm.

\begin{enumerate}
\item For $\LE$-diagram $D$, pick uppermost NW-corner of $1$-rkline, $\alpha=\Delta_{a,b}$.
\item Let $t$ be the maximal rank of the dots in $D$. Weakly-connected property tells us  SW(NE)-end of $(t-1)/\alpha$ and $t/\alpha$-rklines are connected. Pick NE(SW)-end of $t/\alpha$-rkline. Let it be $\Delta_{aX,bY}$.
\item Add $\Delta_{aX,bY}$ to $\S'$. Delete $\alpha$ from $D$ to get $D'$. Repeat the process for $D'$ if it is nonempty. If $D'$ is empty, end by adding $\Delta_{\phi,\phi}$ to $\S'$.
\end{enumerate}

Let's apply this to the first example of Section 7. We show the figures for first few steps.

\begin{figure}[ht]
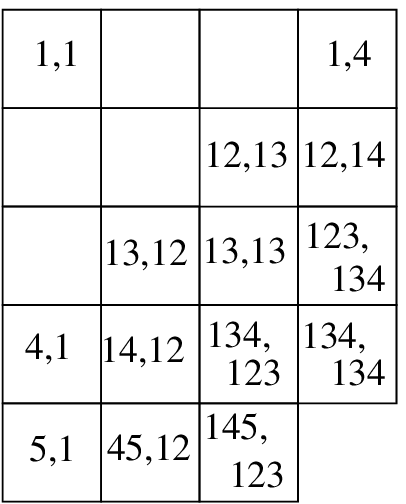
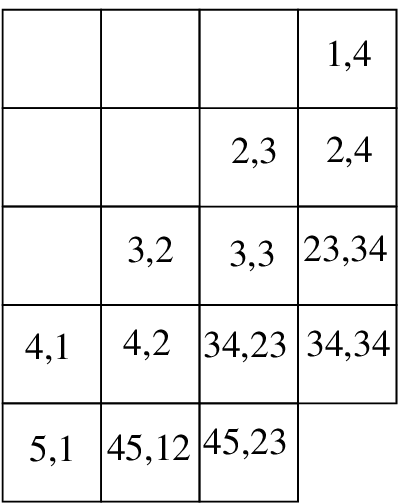
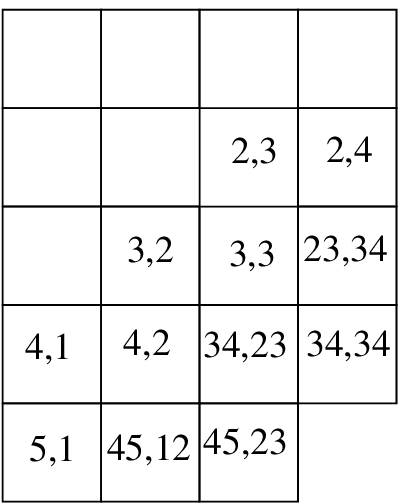
\caption{$\S'=\{\}$, $\S'=\{\Delta_{145,123}\}$, $\S' = \{\Delta_{145,123}, \Delta_{1,4}\}$}
\label{fig:algex-1}
\end{figure}

\begin{figure}[ht]
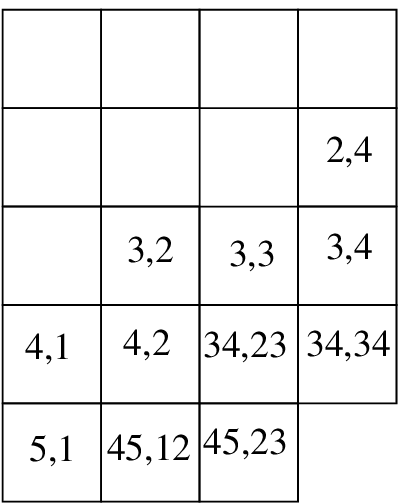
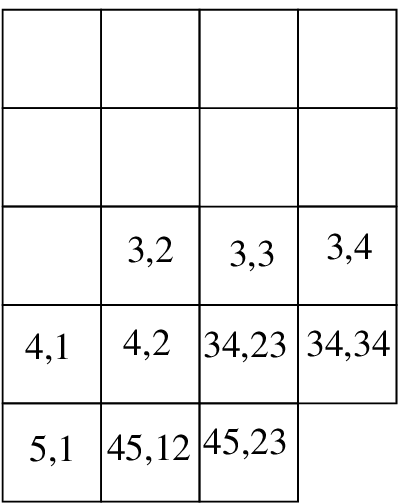
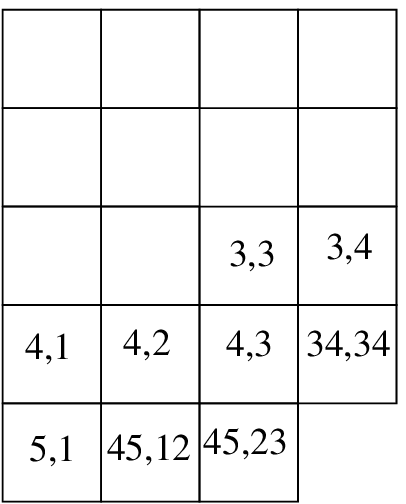
\caption{$\S' = \{\Delta_{145,123}, \Delta_{1,4},\Delta_{23,34}\}$ , $\S' = \{\Delta_{145,123}, \Delta_{1,4},\Delta_{23,34},\Delta_{2,4}  \}$ , $\S' = \{\Delta_{145,123}, \Delta_{1,4},\Delta_{23,34},\Delta_{2,4},\Delta_{34,23}  \}$}
\label{fig:algex-2}
\end{figure}

In the end we would get 
$$\S'  = \{\Delta_{145,123}, \Delta_{1,4},\Delta_{23,34},\Delta_{2,4},\Delta_{34,23},\Delta_{34,34},\Delta_{3,4},\Delta_{45,12},\Delta_{45,23},\Delta_{4,3},\Delta_{4,4},\Delta_{5,1},\Delta_{5,2},\Delta_{5,3},\Delta_{\phi,\phi}  \}$$

It is easy to check that they are all unique-path-variables in $\LE$-network of the cell.

\section{Lattice-path-matroid cells, Cluster algebra and the nonnegativity conjecture}

\begin{figure}[ht]
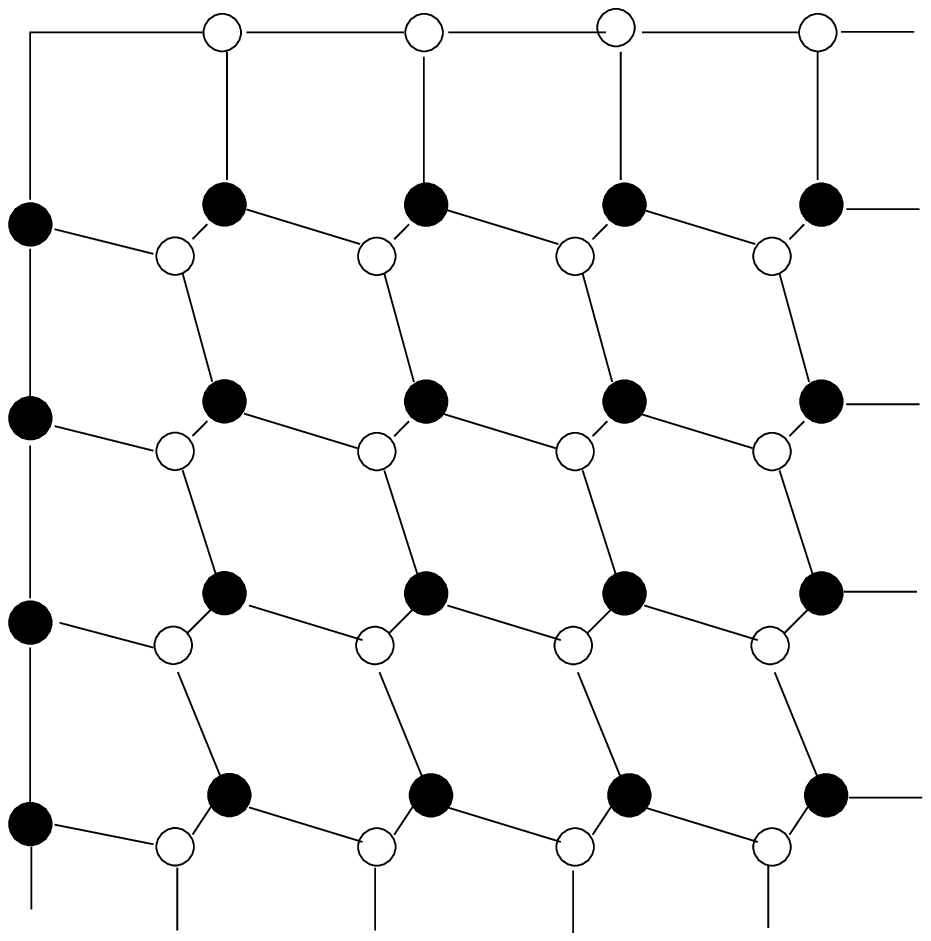
\caption{Plabic graph of top-cell in $Gr_{4,8}^{tnn}$}
\label{fig:trasf-1}

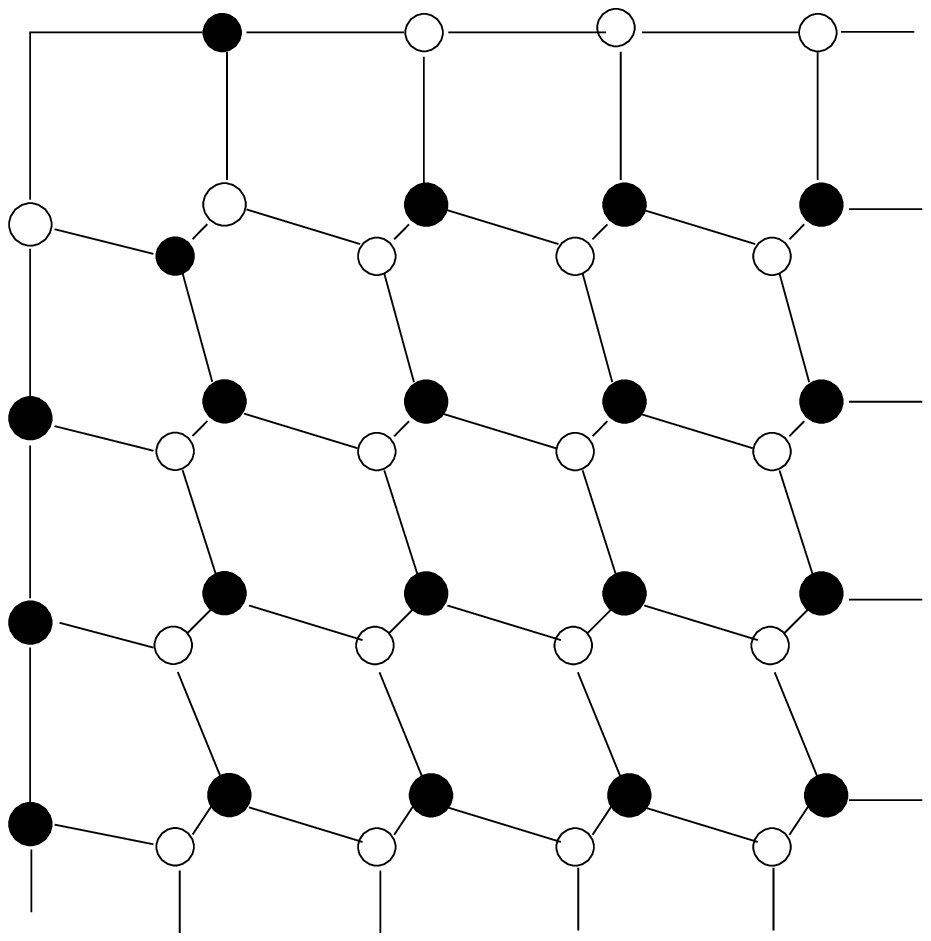
\caption{(M1) on the upper left square}
\label{fig:trasf-2}

\end{figure}

In this section, we will restrict ourselves to lattice-path-matroid cells. We can apply our algorithm for weakly-connected cells to obtain a TP-bases $\S'$ from the canonical TP-bases $\S$. We will prove that $\S'$ is a cluster in the cluster algebra associated to the cell. This will imply that for this seed, every other Pl\"ucker variables are expressed in subtraction-free Laurent polynomials in variables of $\S'$.

Let's see what happens if we apply the algorithm for TP-bases. Pick a NW-corner of $1$-rkline $\alpha$. Then each $k/\alpha$-rkline consists of unique element $\alpha_k$. So we don't need any straightening or refining process. We can jump from $\alpha=\alpha_1$ to $\alpha_2$, from $\alpha_2$ to $\alpha_3$ and so on. Now let's show that this jumping corresponds to seed-mutations. To do this, we want to show this jumping corresponds to 
square moves in a plabic graph, which will prove that the move is a valid seed mutation in the cluster algebra associated to the cell.

So if we restate the algorithm in terms of cluster algebra and quivers, it is as the following. Let the initial seed be $(Q,u)$, pick a northwest corner dot in the $\LE$-diagram. Delete that dot to get a new $\LE$-diagram. And let $(Q',u')$ be the initial seed corresponding to that cell corresponding to the new $\LE$-diagram. Using square moves, we will transform $(Q,u)$ to $(Q'',u'')$ so that it looks almost like $(Q',u')$. That is, excluding one vertice and all edges incident to it in $Q''$, it will look like $Q'$. And $u''$ will consist of variables in $u'$ and one extra variable that corresponds to a unique path in the original $\LE$-graph.

The canonical TP-bases $\S$ (The cluster of inital seed $(Q,u)$) consists of the following Pl\"ucker variables. For each dot at $(a,b)$, denote $t_{a,b}$ as the maximal number such that there is a dot inside the diagram at $(a-t_{a,b},b-t_{a,b})$. Define $E_{a,b}$ as $\{a-t_{a,b},\cdots,a-1,a \}$ and $F_{a,b}$ as $\{b-t_{a,b}, \cdots, b-1,b \}$. Then $\S$ consists of Pl\"ucker variables $\Delta_{E_{a,b},F_{a,b}}$ for each dots in the $\LE$-diagram. Now let's transform this to $\S'$ using square moves of the plabic graph.

First apply the square move to square corresponding to point $(a,b)$. See Figure~\ref{fig:trasf-1} for an example. It is adjacent to squares corresponding to $\Delta_{\phi,\phi},\Delta_{\{a\},\{b+1\}},\Delta_{\{a+1\},\{b\}},\Delta_{\{a,a+1\},\{b,b+1\}}$. So from 
$$\Delta_{\{a\},\{b\}} \Delta_{\{a+1\},\{b+1\}} = \Delta_{\phi,\phi}\Delta_{\{a,a+1\},\{b,b+1\}} + \Delta_{\{a\},\{b+1\}} \Delta_{\{a+1\},\{b\}} $$
we get seed obtained by exchanging $\Delta_{\{a\},\{b\}}$ with $\Delta_{\{a+1\},\{b+1\}}$.

\begin{figure}[ht]
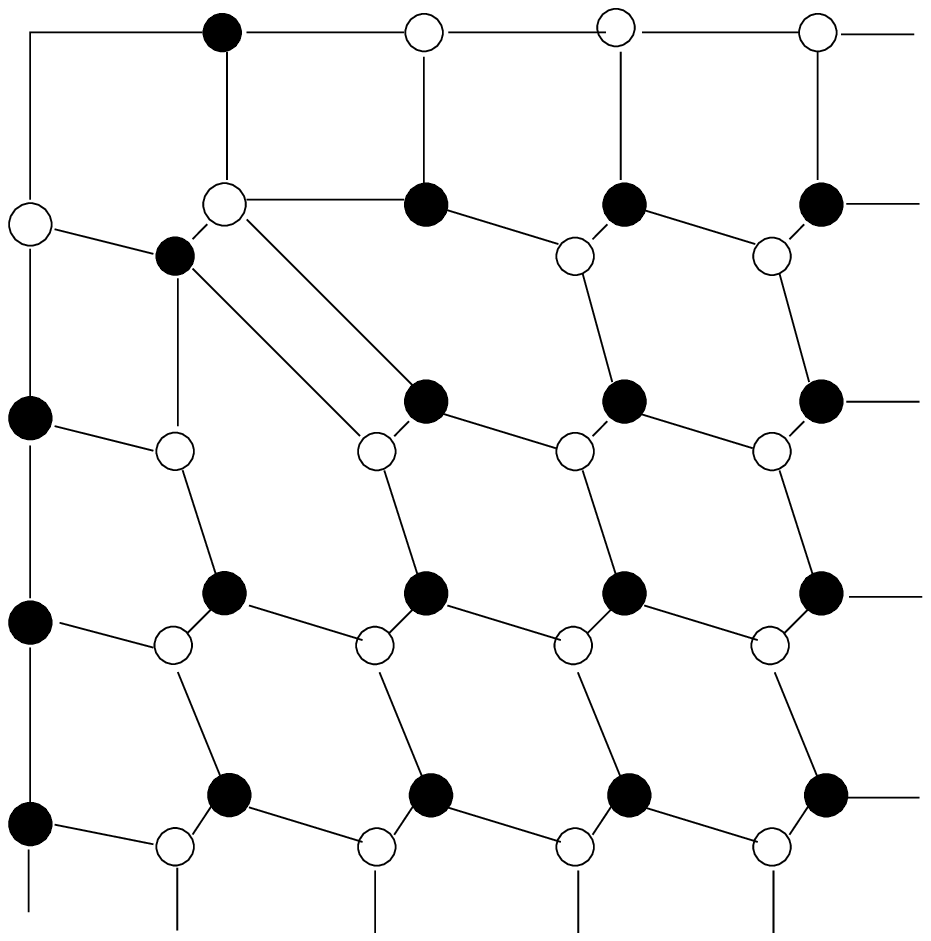
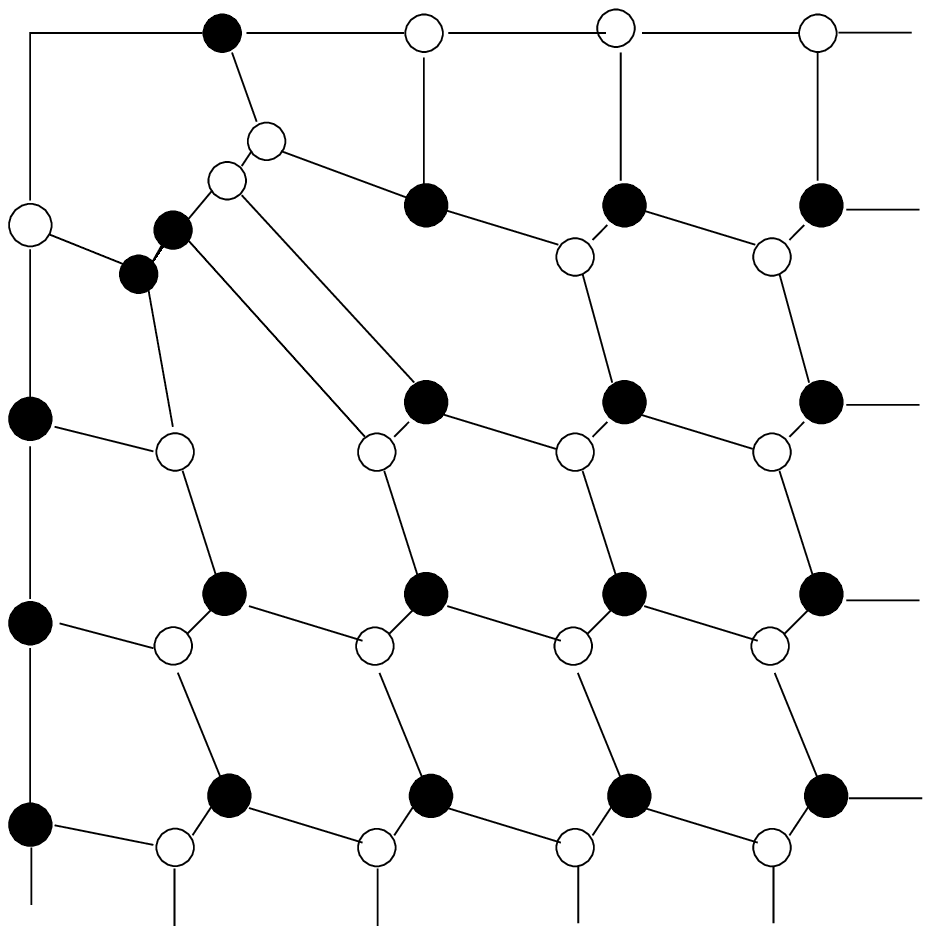
\caption{(M2) to contract two vertices into one. Then (M2) to split a vertice into two vertices.}
\label{fig:trasf-4}
\label{fig:trasf-7}
\end{figure}

\begin{figure}[ht]
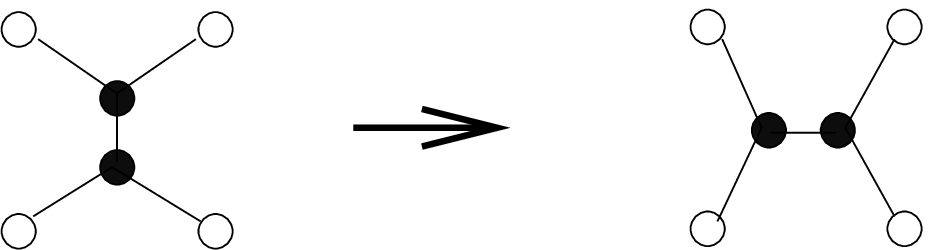
\caption{Contracting then uncontracting using (M2)}
\label{fig:contraf}
\end{figure}

Now for the changed square, pick the black vertice that is on the southeast part. It is connected another black vertice. Contract those two using $(M2)$. Do the same for the white vertice. See Figure~\ref{fig:trasf-2}. Now uncontract the black vertice, but so that the pairs of white vertices that were connected are different from that to prior the contraction. This process is drawn in Figure~\ref{fig:contraf}. Also see Figure~\ref{fig:trasf-4}.

\begin{figure}[ht]
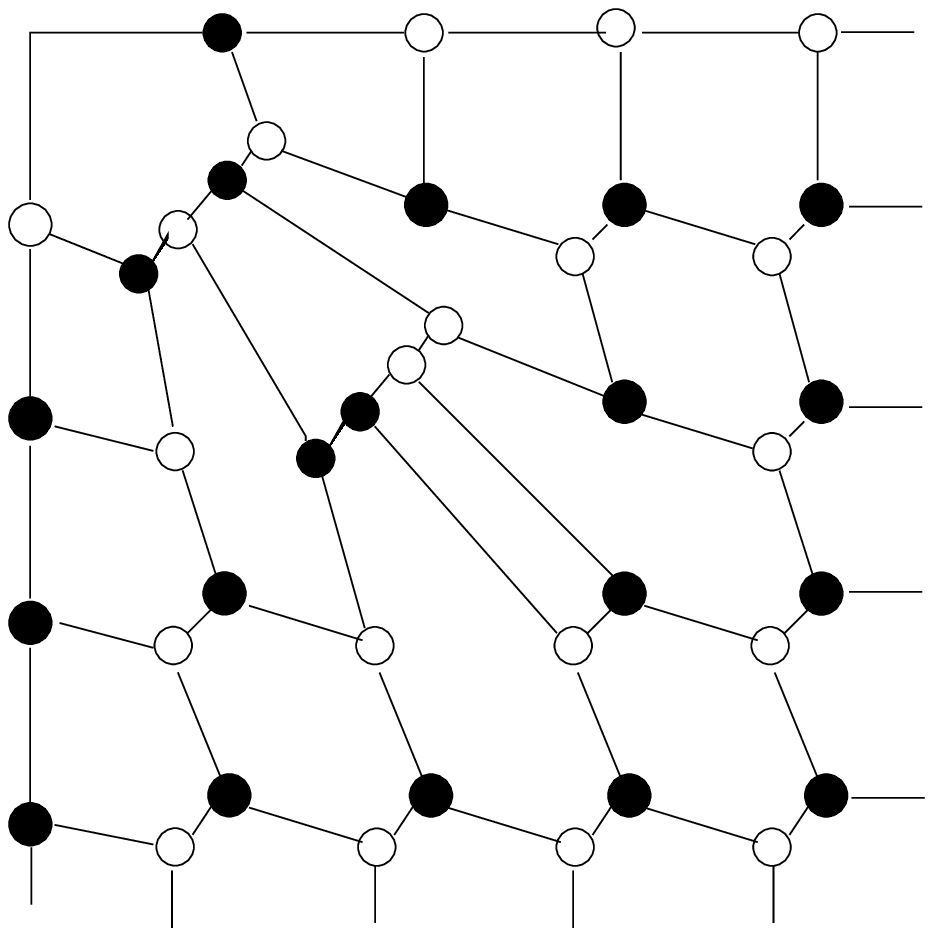
\caption{Similar procedure on the newly formed square, on the lower diagonal of the converted square}
\label{fig:trasf-8}
\end{figure}

Now we can apply the square move to square corresponding to point $(a+1,b+1)$. It is adjacent to squares corresponding to $\Delta_{\{a+1\},\{b+1\}},\Delta_{\{a,a+1\},\{b+1,b+2\}},\Delta_{\{a+1,a+2\},\{b,b+1\}},\Delta_{\{a,a+1,a+2\},\{b,b+1,b+2\}}$. So from 
$$\Delta_{\{a,a+1\},\{b,b+1\}} \Delta_{\{a+1,a+2\},\{b+1,b+2\}} = \Delta_{\{a+1\},\{b+1\}}\Delta_{\{a,a+1,a+2\},\{b,b+1,b+2\}} +  \Delta_{\{a,a+1\},\{b+1,b+2\}} \Delta_{\{a+1,a+2\},\{b,b+1\}} $$
we get seed obtained by exchanging $\Delta_{\{a,a+1\},\{b,b+1\}}$ with $\Delta_{\{a+1,a+2\},\{b+1,b+2\}}$. See Figure~\ref{fig:trasf-8}.

\begin{figure}[ht]
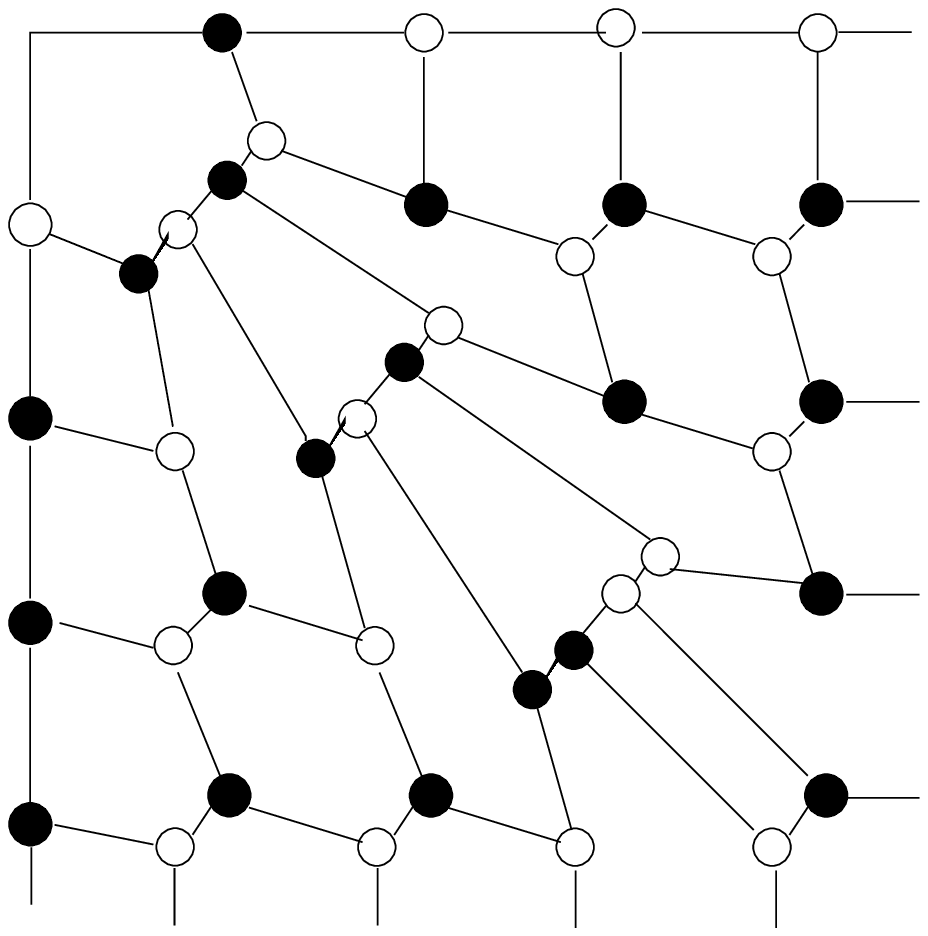
\caption{Repeat the procedure on the newly formed square, end of process}
\label{fig:trasf-9}

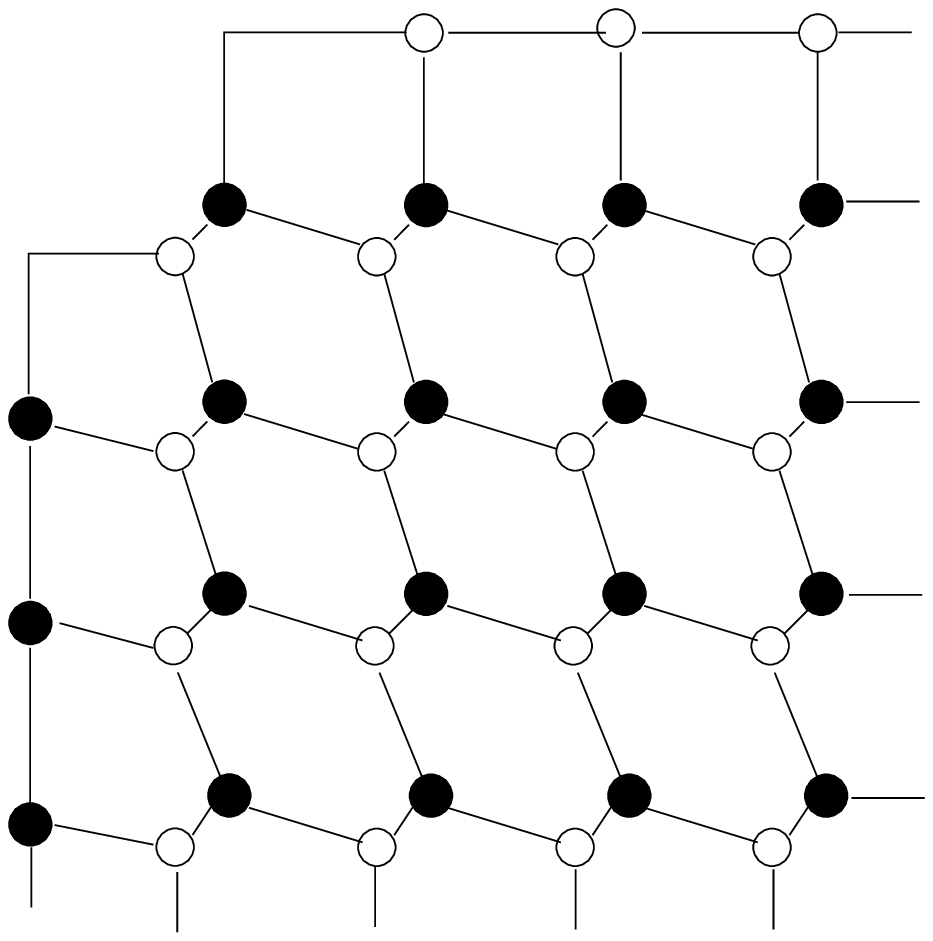
\caption{Plabic graph of cell corresponding to $\LE$-diagram obtained by deleting upper-left dot from $\LE$-diagram of $Gr_{4,8}^{tnn}$}
\label{fig:trasfred}

\end{figure}

Now repeat the procedure, at each step replacing $\Delta_{\{a,a+1,\cdots,a+t\},\{b,b+1,\cdots,b+t\}}$ by \\ $\Delta_{\{a+1,\cdots,a+t+1\},\{b+1,\cdots,b+t+1\}}$ for $t < k$. What we get in the end is the plabic graph of $D_2$ with an extra face attached to the southeastern boundary. See Figure~\ref{fig:trasf-9} and Figure~\ref{fig:trasfred}.

So for each NW-corner of $1$-rkline of the $\LE$-diagram, apply the above procedure and delete that dot. If we finish this procedure for all dots of the $\LE$-diagram, we have obtained $\S'$ from $\S$ using seed mutations. So we have proven the following.

\begin{theorem}
Pick any cell corresponding to lattice-path-matroid in the nonnegative Grassmannian. Look at its $\LE$-diagram. For each box $(a,b)$, denote $k_{a,b}$ as the maximal number such that there is a box inside the diagram at $(a+k_{a,b},b+k_{a,b})$. Define $X_{a,b}$ as $\{a,a+1,\cdots,a+k_{a,b} \}$ and $Y_{a,b}$ as $\{b, b+1, \cdots, b+k_{a,b} \}$.

Let $\S'$ be set of Pl\"ucker variables of from $\Delta_{X_{a,b},Y_{a,b}}$ for all dots in the $\LE$-diagram. Set $x_e$ where $e$ is an edge heading towards $(x,y)$ as $\Delta_{X_x,Y_y}\Delta_{X_{x+1},Y_{y+2}} / \Delta_{X_x,Y_{y+1}}\Delta_{X_{x+1},Y_{y+1}}$.

Since we are looking at projective coordinates, we can set $\Delta_{\phi,\phi}=1$. Then $\Delta_{H,J} = \sum_{P: H \rightarrow J} \prod_{e \in P} x_e$ where $P: H \rightarrow J$ indicates non-vertice-crossing paths from source $H$ to sinks $J$. $\S'$ is a cluster inside the cluster algebra associated to the cell.

\end{theorem}

\begin{proof}
Only thing we haven't proven is to express each $x_e$ using variables of $\S'$. But this is easy and we leave it to the reader as an exercise.
\end{proof}

\section{Further Remarks}

In \cite{DKK}, the tropicalized TP-bases is studied.

\begin{definition}
Let $B$ be a truncated box $B_m^{m'}(a)$, which is defined as set of integer vectors $x=(x_1,\cdots,x_n)$ such that $m \leq x_i \leq max(a_i,m')$ for all $i \in [n]$, $a=(a_1,\cdots,a_n)$. 

A function $f: B \rightarrow \R$ is called a TP-function if it satisfies
$$f(x+1_i+1_k)+f(x+1_j) = max\{ f(x+1_i+1_j)+f(x+1_k),f(x+1_i)+f(x+1_j+1_k) \}$$
for any $x$ and $1 \leq i < j < k \leq n$ and satisfies
$$f(x+1_i+1_k) = f(x+1_j+1_l) = max\{f(x+1_i+1_j)+f(x+1_k+1_l),f(x+1_i+1_l)+f(x+1_j+1_k)\}$$
for any $x$ and $1 \leq i < j < k < l \leq n$ provided all six vectors occurring in each relation belong to $B$. 

Then a subset $\B \subseteq B$ is called a TP-basis if each TP-function on B is determined by values on $\B$, and values on $\B$ can be chosen arbitrarily.
\end{definition}
 When we set $m=0,m'=1,a=(1,\cdots,1)$ this corresponds to the set of Pl\"ucker variables in the top-cell of the nonnegative Grassmannian, since the above relations are tropical analogues of the 3-term Pl\"ucker relation. The TP-bases we used, is different in the sense that instead of looking at functions like TP-functions, we restrict ourselves only to Pl\"ucker variables. But we look not only at the top cell, but all other positroid cells in the nonnegative Grassmannian.

\medskip

\begin{figure}[ht]
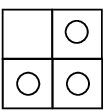
\caption{An example of a cell when set of nonzero Pl\"ucker variables $\not \subseteq$ set of cluster variables}
\label{fig:weird}
\end{figure}

Now let's try comparing set of nonzero Pl\"ucker variables of the cell and set of cluster variables of the associated cluster algebra. In general, the set of nonzero Pl\"ucker variables of the cell are not contained in the set of cluster variables of the associated cluster algebra. An easy example of such cell is a cell having $\LE$-diagram in Figure~\ref{fig:weird}. The set of cluster variables in this cell is $\{\Delta_{2,2},\Delta_{1,2},\Delta_{2,1},\Delta_{\phi,\phi}\}$. But the set of nonzero Pl\"ucker variables in this cell is $\{\Delta_{2,2},\Delta_{1,2},\Delta_{2,1},\Delta_{1,1},\Delta_{\phi,\phi}\}$.

\begin{conjecture}
\label{conj:allcell}
For any positroid cell, the canonical TP-bases can be transformed via Pl\"ucker mutation to a TP-bases consisting of variables coming from unique paths.
\end{conjecture}

The conjecture has been checked for all cells inside $Gr_{kn}^{tnn}, n \leq 8$. We have proven the conjecture for weakly-connected cells using Lemma~\ref{lem:main}. For non weakly-connected cells we can't use the approach of deleting dots one-by-one as we did to prove Theorem~\ref{thm:main}. So we need a different approach.

\begin{conjecture}
For any positroid cell, all TP-bases are mutation equivalent.
\end{conjecture}

If the above conjecture is true, then the following conjecture may also be true.

\begin{conjecture}
For any cell, pick any TP-bases. Then all other Pl\"ucker variables are expressible as subtraction-free Laurent polynomials using the variables in that TP-bases.
\end{conjecture}

This conjecture is the analogue of the nonnegativity conjecture of cluster algebras. 
Similar conjecture was given in \cite{FZ6} for TP-bases coming from double wiring diagrams, inside cells of nonnegative part of $GL_k$.  If we prove Conjecture~\ref{conj:allcell}, we have a TP-bases such that all other nonzero Pl\"ucker variables can be expressed as subtraction-free Laurent polynomials in variables of the TP-bases. If we have used mutations of form $\Delta_{H_1}\Delta_{H_2} = \Delta_{H_3}\Delta_{H_4}$, then the new set of TP-bases also has the property. So we would have to prove that whenever a TP-bases $S_1$ has the property and we obtained $S_2$ from $S_1$ by $\Delta_{H_1}\Delta_{H_2} = \Delta_{H_3}\Delta_{H_4} + \Delta_{H_5}\Delta_{H_6}$, $S_2$ also has the property.


\begin{thebibliography}{HLSS}

\bibitem[BGW]{BGW} A.~Borovik, I.~Gelfand, N.~White. \textit{Coxeter Matroids} Birkhauser, Boston, 2003.

\bibitem[BLSWZ]{BLSWZ} A.~Bj\"orner, M. Las Vergnas, B. Sturmfels, N. White, and  G.~M.~Ziegler.
\textit{Oriented Matroids}. Encyclopedia of Mathematics and Its Applications,
vol. 46. Cambridge University Press, Cambridge, 1993.

\bibitem[BMN]{BMN} J.~Bonin, A.~de Mier, M.~Noy. Lattice path matroids: enumerative aspects and Tutte polynomials,\textit{ J. Combin. Theory Ser. A} 104 (2003) 63-94.

\bibitem[CR]{CR} P. Caldero, M. Reineke. On the quiver Grassmannian in the acyclic case, 
arXiv: math/0611074 [math.RT].

\bibitem[CP]{CP} G. Carroll and G. Price, Two new combinatorial models for the Ptolemy recurrence, unpublished
memo (2003)

\bibitem[DKK]{DKK} V.~Danilov, A.~Karzanov, G.~Koshevoy. On bases of tropical Pl\"ucker functions, arXiv:0712.3996v2 [math.CO]

\bibitem[F]{F} W.~Fulton. \textit{Young Tableaux. With Applications to Representation Theory and Geometry} New York: Cambridge University Press, 1997

\bibitem[FZ]{FZ} S.~Fomin, A.~Zelevinsky, Cluster algebras. I. Foundations, \textit{J. Amer. Math. Soc}. 15 (2002), 497-529.

\bibitem[FZ2]{FZ2} S.~Fomin, A.~Zelevinsky, Cluster algebras. II. Finite type classification, \textit{Invent. Math}. 154 (2003), 63-121.

\bibitem[FZ3]{FZ3} S.~Fomin, A.~Zelevinsky, The Laurent phenomenon, \textit{Adv. Applied Math}. 28 (2002), 119-144.

\bibitem[FZ4]{FZ4} S.~Fomin, A.~Zelevinsky, Cluster algebras: notes for the CDM-03 conference, \textit{Current developments in mathematics}, 2003, Int. Press, Somerville, MA, (2003), 1-34.

\bibitem[FZ5]{FZ5} S.~Fomin, A.~Zelevinsky, Double Bruhat cells and total positivity, \textit{J. Amer. Math. Soc}. 12 (1999), 335-380.

\bibitem[FZ6]{FZ6} S.~Fomin, A.~Zelevinsky, Total positivity: Tests and parameterizations, \textit{The Mathematical Intelligencer}. 22 (2000), 23-33.


\bibitem[G]{G} D.~Gale, Optimal assignments in an ordered set: an application of matroid theory, \textit{J. Combinatorial Theory} \textbf{4}(1968), 1073-1082.

\bibitem[HS]{HS} A.~Henriques, D.~Speyer, The multidimensional cube recurrence, arXiv: math/0708.2478 [math.CO].

\bibitem[K]{K} B.~Keller. Cluster algebras, quiver representations and triangulated categories, arXiv:0807.1960v2 [math.RT].


\bibitem[LI]{LI} B.~Lindstrom. On the vector representations of induced matroids, \textit{Bull. London Math. Soc}. 5 (1973), 85-90.

\bibitem[M]{M} G.~Musiker, A graph theoretic expansion formula for cluster algebras of classical type, arXiv: math/0710.3574 [math.CO].

\bibitem[O]{O} S.~Oh. Positroids and Schubert matroids, arXiv:0803.1018 [math.CO].

\bibitem[P]{P} A.~Postnikov. Total positivity, Grassmannians, and networks, 
arXiv: math/ 0609764v1 [math.CO].
 
\bibitem[PW]{PW} A.~Postnikov, L.~Williams. \textit{Private communication}.

\end{thebibliography}
\end{document}